\documentclass{amsart}
\usepackage[utf8]{inputenc}
\usepackage{amsmath}
\usepackage{amssymb}
\usepackage{caption}
\usepackage{amsthm}
\usepackage[backref=false, colorlinks, linkcolor=blue, citecolor=blue]{hyperref}
\usepackage{cleveref}
\usepackage{hyperref}
\usepackage{nicefrac}
\usepackage[inline]{enumitem}
\usepackage{enumitem}

\usepackage{changepage} 

\usepackage{todonotes}
\usepackage{amssymb}
\usepackage{mathrsfs}

\usepackage{array,float}

\usepackage{tikz}
\usetikzlibrary{
  knots,
  hobby,
  decorations.pathreplacing,
  shapes.geometric,
  calc
}
\usetikzlibrary{shapes,arrows}
\usetikzlibrary{fit,positioning}
\usetikzlibrary{patterns,decorations.pathreplacing}
\usepackage{xkeyval}
\usepackage{moreverb}
\usepackage{epic}
\usepgfmodule{shapes,plot,decorations}

\usepackage{booktabs} 
\usepackage[normalem]{ulem}

\newcommand{\Z}{\mathbb{Z}}
\newcommand{\Q}{\mathbb{Q}}
\newcommand{\R}{\mathbb{R}}
\newcommand{\CC}{\mathbb{C}}

\newcommand{\HH}{\mathbb{H}}

\newcommand{\PO}{\mathbf{PO}}
\newcommand{\G}{\mathbf{G}}

\usepackage[warn]{mathtext}

\newcommand{\D}{\mathcal{D}}

\newcommand{\E}{\mathbb{E}}

\newcommand{\OOO}{\mathcal{O}}

\newcommand{\id}{\mathrm{Id}}

\newcommand{\Isom}{\mathrm{Isom}}
\newcommand{\Cyc}{\mathrm{Cyc}}

\newcommand{\Comm}{\mathrm{Comm}}

\newcommand{\tr}{\mathrm{tr}}

\newtheorem{theorem}{Theorem}[section]
\newtheorem{corollary}{Corollary}[theorem]
\newtheorem{lemma}[theorem]{Lemma}
\newtheorem{definition}[theorem]{Definition}

\newtheorem{question}[theorem]{Question}

\theoremstyle{remark}
\newtheorem{remark}[theorem]{Remark}

\theoremstyle{remark}

\renewcommand{\qed}{\hfill$\scriptstyle\blacksquare$}

\makeatletter
\def\@tocline#1#2#3#4#5#6#7{\relax
  \ifnum #1>\c@tocdepth 
  \else
    \par \addpenalty\@secpenalty\addvspace{#2}%
    \begingroup \hyphenpenalty\@M
    \@ifempty{#4}{%
      \@tempdima\csname r@tocindent\number#1\endcsname\relax
    }{%
      \@tempdima#4\relax
    }%
    \parindent\z@ \leftskip#3\relax \advance\leftskip\@tempdima\relax
    \rightskip\@pnumwidth plus4em \parfillskip-\@pnumwidth
    #5\leavevmode\hskip-\@tempdima
      \ifcase #1
       \or\or \hskip 1em \or \hskip 2em \else \hskip 3em \fi%
      #6\nobreak\relax
    \dotfill\hbox to\@pnumwidth{\@tocpagenum{#7}}\par
    \nobreak
    \endgroup
  \fi}
\makeatother

\title[On quasi-arithmeticity of hyperbolic gluings]{On quasi-arithmeticity of hyperbolic gluings}

\author{Nikolay Bogachev}
\address{Department of Computer and Mathematical Sciences, University of Toronto Scarborough, 1095 Military Trail, Toronto, ON M1C 1A3, Canada}
\email{n.bogachev@utoronto.ca}

\author{Dmitry Guschin}
\address{Moscow Institute of Physics and Technology, Moscow, Russia}
\email{dmckg1999@gmail.com}

\author{Andrei Vesnin}
\address{Sobolev Institute of Mathematics, Novosibirsk, Russia}
\address{Tomsk State University, Tomsk, Russia}
\email{vesnin@math.nsc.ru }

\begin{document}

\begin{abstract}
We study a more general version of the gluings of hyperbolic orbifolds in the spirit of Gromov and Piatetski-Shapiro, where the gluing pieces, called the building blocks, are no longer assumed to be arithmetic or incommensurable. 

We prove that if such a general hyperbolic gluing along a common finite-volume totally geodesic hypersurface is quasi-arithmetic (this is a broader notion than that of arithmeticity) then each building block must be quasi-arithmetic as well and, moreover, with the same ambient group and adjoint trace field. We also show that there exist arithmetic gluings whose building blocks are incommensurable even despite the reflection with respect to the lift of the gluing locus commensurates the fundamental group of the gluing. On the other hand, we provide an example of nonarithmetic but quasi-arithmetic orbifolds such that a specific gluing of such an orbifold with itself along the boundary gives rise to an arithmetic hyperbolic orbifold.

We illustrate the above results in the setting of reflection groups and hyperbolic Coxeter polyhedra and apply them to rule out the (quasi-)arithmeticity of a family of ideal hyperbolic right-angled $3$-polyhedra, namely, certain ``twisted'' ideal right-angled antiprisms, which play an important role in low-dimensional geometry and topology. 
\end{abstract}


\maketitle


\section{Introduction}

In recent years, there has been growing interest in the study of different arithmetic types of hyperbolic lattices,  in large part due to the ability of the arithmetic properties to distinguish commensurability classes of such lattices and of the associated hyperbolic manifolds and orbifolds. We are particularly interested in the class of {\em quasi-arithmetic} hyperbolic lattices, as well as in the so-called {\em hyperbolic gluings}, i.e. the hyperbolic orbifolds obtained by gluing some number of finite-volume hyperbolic orbifolds along their isometric boundaries. The study of arithmetic properties of hyperbolic gluings is motivated by the ideas of Gromov and Piatetski-Shapiro \cite{GPS87} for constructing nonarithmetic hyperbolic lattices (including ones generated by reflections, cf. Vinberg \cite{Vin14}).

The notion of {\em quasi-arithmeticity}, introduced by Vinberg in 1967 in his seminal work \cite{Vin67} on hyperbolic reflection groups, is indeed broader than that of arithmeticity. The first quasi-arithmetic but not arithmetic examples were provided by Vinberg in low dimensions in the same paper \cite{Vin67}, but nowadays such examples are known to exist in all dimensions. Thomson \cite{Tho16} showed that the nonarithmetic manifolds constructed by Agol~\cite{agol2006systoles}, Belolipetsky--Thomson \cite{BT11}, and Bergeron--Haglund--Wise \cite{BHW11} (see also Mila~\cite{Mila18}), giving rise to finite-volume hyperbolic manifolds with arbitrarily small systole, are quasi-arithmetic. On the other hand, it was demonstrated in the same work \cite{Tho16} of Thomson that the initial gluing construction of Gromov and Piatetski-Shapiro \cite{GPS87} gives rise to non-quasi-arithmetic lattices (commonly called {\em classical hybrids}, cf. Mila \cite{Mila-imrn}). 

Our main result, Theorem~\ref{th:general-GPS}, concerns the (quasi-)arithmetic properties of hyperbolic gluings and their building blocks.

Let $\Gamma_1,\Gamma_2 < \Isom(\mathbb{H}^n)$ be two lattices both containing a reflection $r$ in a hyperplane $H \subset \mathbb{H}^n$ such that:

\begin{enumerate}
    \item for any $\gamma \in \Gamma_i$, $i=1,2$, either $\gamma(H)=H$ or $\gamma(H) \cap H = \emptyset$,
    \item $N_{\Gamma_1}(r)=N_{\Gamma_2}(r)= \langle r \rangle \times \Gamma_0$, where $\Gamma_0$ leaves invariant the two half--spaces bounded by $H$ and acts as a lattice in $H$. Here, $N_{\Gamma_i}(r)$ denotes the normalizer of $r$ in $\Gamma_i$.    
\end{enumerate} 

This means that the projections of $H$ in finite-volume hyperbolic orbifolds $\HH^n/\Gamma_1$ and $\HH^n/\Gamma_2$ are isometric embedded finite-volume totally geodesic subspaces. 

The following is an explicit group-theoretic description of the initial hybrid construction of Gromov and Piatetski--Shapiro \cite{GPS87}, see also Vinberg \cite{Vin14}.
Consider the set $\Omega_i = \{\gamma(H) \mid \gamma \in \Gamma_i\}$ for a fixed $i = 1$ or $2$. The set $\Omega_i$ decomposes the space $\mathbb{H}^n$ into a collection of closed domains transitively permuted by the group $\Gamma_i$, and each of these pieces is a fundamental domain for the action of a reflection group $N_i = \langle \gamma r \gamma^{-1} \mid \gamma \in \Gamma_i \rangle$, which is the normal closure of $r$ in $\Gamma_i$. Let $D_i$ be one of these pieces, and let $\Delta_i=\{\gamma \in \Gamma_i \mid \gamma(D_i)=D_i\}$. Then the group $\Gamma_i$ can be decomposed into a semidirect product $\Gamma_i=N_i\rtimes \Delta_i$. Moreover, by choosing $D_1$ and $D_2$ to belong to opposite half-spaces with respect to the hyperplane $H$, we can ensure that $\Delta_1 \cap \Delta_2=\Gamma_{0}$ is the common stabilizer of the hyperplane $H$ in both lattices $\Gamma_1$ and $\Gamma_2$.

Theorem~\ref{th:general-GPS} provides the sufficient conditions for building a new {\em gluing} non-quasi-arithmetic lattice out of the lattices $\Gamma_1$ and $\Gamma_2$ described above. To formulate this theorem we recall that Vinberg \cite{Vin71} introduced two commensurability invariants for lattices in semi-simple algebraic groups: the {\em adjoint trace field} $k_\Gamma$ and the {\em ambient $k$-group} $\G_\Gamma$ (see Section~\ref{sec:trace-fields} for details).

\begin{theorem}\label{th:general-GPS}
Given $\Gamma_i$, $i=0,1,2$, and $\Delta_j$, $j=1,2$, as above, the following conditions hold:
\begin{itemize}
    \item[\textnormal{(1)}]The group $\Gamma = \langle \Delta_1, \Delta_2 \rangle = \Delta_1 *_{\Gamma_0} \Delta_2 < \mathrm{Isom}(\mathbb{H}^n)$ is a lattice (it will be referred to as a gluing lattice).
    \item[\textnormal{(2)}] If $\Gamma$ is quasi-arithmetic (in particular, if it is arithmetic) in $\mathrm{Isom}(\mathbb{H}^n)$ with ambient group $\G$ and adjoint trace field $k$, then both $\Gamma_1$ and $\Gamma_2$ are quasi-arithmetic with the same $\G$ and $k$.
    \item[\textnormal{(3)}] If the lattice $\Gamma$ is nonarithmetic and $r \in \Comm(\Gamma)$, then $\Gamma_1$ and $\Gamma_2$ are commensurable. 
\end{itemize}
\end{theorem}

\noindent Part (2) of Theorem \ref{th:general-GPS} has the following obvious but useful technical consequence.
\begin{corollary}\label{cor:non-quasi-fields}
Let $\Gamma = \langle \Delta_1, \Delta_2 \rangle = \Delta_1 *_{\Gamma_0} \Delta_2 < \mathrm{Isom}(\mathbb{H}^n)$ be a gluing lattice defined as above. If
\begin{itemize}
    \item[\textnormal{(1)}] the adjoint trace fields $k_1$ and $k_2$ of $\Gamma_1$ and $\Gamma_2$, respectively, do not coincide, or
    \item[\textnormal{(2)}] $k_1 = k_2 = k$, but ambient groups $\mathbf{G}_1$ and $\mathbf{G}_2$ of $\Gamma_1$ and $\Gamma_2$, respectively, are not $k$-isomorphic,
\end{itemize}
then $\Gamma$ is not quasi-arithmetic.
\end{corollary}

\begin{remark}\label{rem:gluing}
In the terminology of geometric topology, the construction from Theorem \ref{th:general-GPS} means that we {\em glue} a finite-volume hyperbolic orbifold $M=\HH^n/\Gamma$ from two finite-volume hyperbolic orbifolds $M'_1$ and $M'_2$ (called {\em building blocks}) with a common finite-volume totally geodesic boundary $\partial M'_1 = \partial M_2'$, which is fixed pointwise by the reflection $r \in \Gamma_1 \cap \Gamma_2$. Here, each $M'_i$, $i=1,2$, is obtained from the orbifold $M_i = \HH^n/\Gamma_i$ by ``forgetting'' the fact that the totally geodesic finite-volume codimension-$1$ suborbifold $N = H/\Gamma_0$ is the singular locus of the reflection $r \in \Gamma_1 \cap \Gamma_2$. However, we will frequently conflate $M_i$ and $M_i'$ when speaking about building blocks. See Figure~\ref{fig:gluing-proof} for the visualization of $M$ and $M'_i$ in the setting of manifolds.
\end{remark}

The trace fields of gluings of arithmetic pieces were studied by Mila \cite{Mila-imrn}. It is also worth mentioning that Emery--Mila \cite{EM21} proved that such gluings are pseudo-arithmetic (this broader notion than that of quasi-arithmeticity is not in the scope of the present paper).

It turns out that the nonarithmeticity assumption in Theorem~\ref{th:general-GPS}, part (3), is crucial, as there exist  arithmetic lattices $\Gamma$ built out of incommensurable lattices $\Gamma_i$ even when the reflection along the gluing locus commensurates $\Gamma$ (i.e. when $r \in \Comm(\Gamma)$).

\begin{theorem}\label{th:arith-surfaces-incomm-blocks}
    Let $S_{g,1}$ be an arithmetic hyperbolic surface of genus $g \ge2$ with one cusp. Then any separating simple closed geodesic $\alpha$ on $S_{g,1}$ provides the gluing of $S_{g,1}$ from two orbifolds $M_1$ and $M_2$ with boundaries $\partial M_1 \cong \partial M_2 \cong \alpha$ such that the reflection $r$ with respect to the geodesic line $\widetilde{\alpha}$, being the lift of  $\alpha$ in $\HH^2$, satisfies $r \in \Comm(\pi_1 S_{g,1})$, but $\Gamma_1$ and $\Gamma_2$ are incommensurable. 
\end{theorem}

\begin{figure}
   \begin{center}
        \includegraphics[width=0.7\linewidth]{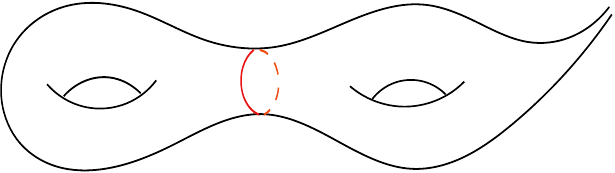}
        \caption{The red arc $\alpha$ splits the surface $S_{2,1}$ into two orbifolds with boundary $\alpha$, which is the singular locus of the reflection $r$ with respect to the lift $\widetilde{\alpha}$ in $\HH^2$. Here we see that the left building block of $S_{2,1}$ is compact, but the right one is not.}
        \label{fig:S2-1}
    \end{center}
    \end{figure}

See Figure \ref{fig:S2-1} for the case of $S_{2,1}$ in Theorem~\ref{th:arith-surfaces-incomm-blocks}.

\medskip

Let us illustrate the gluing construction from Theorem~\ref{th:general-GPS} in the setting of reflection orbifolds and hyperbolic Coxeter polyhedra.

\begin{theorem}\label{th:Coxeter-orbifolds}
Let $P_1$ and $P_2$ be two finite-volume Coxeter polyhedra in $\HH^n$. Suppose that $P_1$ and $P_2$ have a common facet $F$ (or an edge of finite length if $n=2$) which meets all of its adjacent facets both in $P_1$ and $P_2$ at even angles, i.e. angles of the form $\pi/2m$, $m \geqslant 1$, such that the corresponding adjacent (with respect to the common facet $F$) angles are equal for $P_1$ and $P_2$.
 
Then the group $\Gamma_P$ is a lattice, and it may be quasi-arithmetic only in the case when both $\Gamma_{P_1}$ and $\Gamma_{P_2}$ are quasi-arithmetic with the same adjoint trace field and ambient group as $\Gamma_P$ itself.
\end{theorem}

\begin{figure}
    \centering
    \includegraphics[scale=1.15]{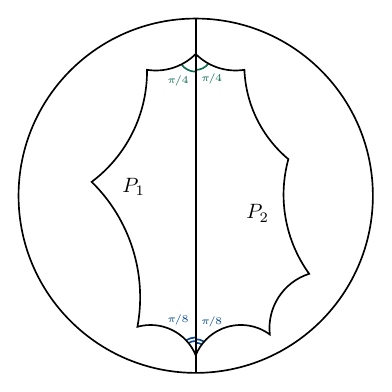}
    \caption{A Coxeter polygon $P$ glued from two other Coxeter polygons $P_1$ and $P_2$ having a common side $F$. This side $F$ has angles $\pi/4$ and $\pi/8$ with adjacent sides in both polygons.}
    \label{fig:hyb-Coxeter}
\end{figure}

\noindent A $2$-dimensional illustration of the angle condition in Theorem~\ref{th:Coxeter-orbifolds} can be found in Figure~\ref{fig:hyb-Coxeter}.

A very special case of the gluing construction is a simple doubling trick that allows one to produce infinitely many right-angled polyhedra by doubling a single right-angled polyhedron along its walls. This method was successfully applied by Allcock \cite{All06} to a more general family of polyhedra than those that are right-angled. We would like to remark that Makarov \cite{Mak68} constructed infinitely many compact Coxeter polyhedra precisely as in Theorem \ref{th:Coxeter-orbifolds}, by gluing two different polyhedra along their common facets. It is worth mentioning that quasi-arithmeticity played an important role in the recent paper \cite{BDR24} showing that Makarov's \cite{Mak68} polyhedra give rise even to infinitely many commensurability classes of cocompact reflection groups in $\HH^4$ and $\HH^5$.  Similar gluing (or in fact hybrid) construction was also used by Ruzmanov \cite{Ruz89} and later by Vinberg \cite{Vin14} to provide nonarithmetic hyperbolic reflection groups in higher dimensions.

The following result shows that one can even construct an arithmetic gluing out of commensurable {\em properly quasi-arithmetic} lattices (i.e. quasi-arithmetic but not arithmetic). 

\begin{theorem}\label{th:arith-out-of-quasi}
There exist nonarithmetic (in fact, properly quasi-arithmetic) finite-volume hyperbolic $3$-orbifolds and $3$-manifolds $M$ with totally geodesic boundary $\partial M$ such that the gluing of $M$ with a certain twist $\tau$ of $M$ along $\partial M$ gives an arithmetic hyperbolic $3$-orbifold $N = M \cup_\tau M$. 

In particular, there exists a properly quasi-arithmetic Coxeter polyhedron $P$ in $\HH^3$ with a facet $F$ orthogonal to all its adjacent facets such that an appropriate gluing of $P$ with itself along $F$ gives rise to an arithmetic Coxeter polyhedron $P'$.
\end{theorem}
\begin{remark}
Obviously, $N$, being arithmetic, is not commensurable with a standard reflective double $N'$ of $M$ along $\partial M$, since the latter double $N'$ is nonarithmetic (but quasi-arithmetic in view of Theorem~\ref{th:general-GPS}) as $M$. Similarly, $P'$ is not commensurable with a standard reflective double of $P$ along $F$.
\end{remark}

Now we apply Theorem~\ref{th:general-GPS} and Theorem~\ref{th:Coxeter-orbifolds} to {\em right-angled reflection groups} and their fundamental polyhedra.
The interest in study of geometric and arithmetic properties of {\em hyperbolic right-angled polyhedra} in large part is based on their connection to hyperbolic $3$-manifolds. There is a well-studied infinite family of compact hyperbolic right-angled polyhedra known as {\em L\"obell polyhedra}, where the simplest instance is a dodecahedron. Arithmeticity of Coxeter groups generated by reflections in faces of the L\"obell polyhedra was studied in \cite{AMR09, Ves91}. It was shown in \cite{BD23} that these polyhedra give rise to particularly straightforward examples of closed hyperbolic $3$-manifolds with arbitrarily small systole, and constitute an infinite family even up to commensurability. 

In 1991, Reid \cite{Reid} proved that the figure-eight knot complement is the only  arithmetic hyperbolic knot complement. On the other hand, the classification of arithmetic link complements is still far from complete. Many hyperbolic links, in particular, the {\em fully augmented links} are associated to {\em ideal right-angled polyhedra} in $\HH^3$, see~\cite{Lack04, Pur11, Pur20}. Fundamental groups of such link complements are not always commensurable with the associated right-angled reflection groups, but in many cases they are, see \cite{CDBW12}. Thus, it is desirable to classify arithmetic ideal hyperbolic right-angled $3$-dimensional polyhedra where we mean the arithmeticity of the corresponding reflection groups. 

For instance, the class of ideal right-angled polyhedra contains an infinite family of {\em antiprisms} $A_n$ with the simplest example being an octahedron $A_3$. Let $\D_{2n}$ denote the classical chain link with $2n$ components. It is known that its complement $M_n = \mathbb{S}^3 \setminus \D_{2n}$ admits a decomposition into four isometric copies of the ideal right-angled antiprism $A_n$ and that $\pi_1 M_n$ is commensurable with the reflection groups associated to $A_n$. Recently, Meyer--Millichap--Trapp \cite{MMT20} and independently Kellerhals \cite{Kel23} showed that the ideal right-angled antiprisms $A_n$ provide a sequence of pairwise incommensurable reflection groups, and proved that the corresponding reflection groups are arithmetic if and only if $n = 3$ or $4$. These results imply that $M_n = \mathbb{S}^3 \setminus \D_{2n}$ is an arithmetic hyperbolic $3$-manifold only for the same values: $n=3$ and $n=4$.

Ideal right-angled antiprisms $A_n$, $n \geqslant 3$, constitute an especially interesting subclass of ideal hyperbolic right-angled polyhedra also due to the fact that any ideal hyperbolic right-angled polyhedron can be obtained from one of the antiprisms by a series of so-called {\em edge twist} operations, see Vesnin \cite[Theorem 2.14]{Ves17} (cf. Erokhovets \cite[Theorem 9.13]{Er19}). This is a combinatorial operation: one takes two disjoint edges belonging to a single face, removes them, then adds a new vertex and joins this vertex with the vertices of the removed edges (see Figure~\ref{fig:edge-twist} where all vertices are four-valent and the edges participating in twisting are red). The resulting combinatorial polyhedron can always be realized in $\HH^3$ as an ideal right-angled polyhedron by Andreev's theorem.

\begin{figure}[ht]
    \centering
    \begin{tikzpicture}[scale=1.2] 
\coordinate (1) at (0, 0) ;
  \coordinate (2) at (1, 0) ;
  \coordinate (3) at (0, 1) ;
  \coordinate (4) at (1, 1) ;
  \draw[line width=2.pt, red] (1) -- (3) ;
  \draw[line width=2.pt, red] (2) -- (4) ; 
    \draw[thick](0,1) -- (0.2,1.2);
    \draw[thick](0,1) -- (-0.2,1.2);
    \draw[thick](0,1) -- (-0.2,0.8);
    \draw[thick](1,1) -- (0.8,1.2);
    \draw[thick](1,1) -- (1.2,1.2);
    \draw[thick](1,1) -- (1.2,0.8);
    \draw[thick](1,0) -- (1.2,0.2);
    \draw[thick](1,0) -- (1.2,-0.2);
    \draw[thick](1,0) -- (0.8,-0.2);
    \draw[thick](0,0) -- (0.2,-0.2);
    \draw[thick](0,0) -- (-0.2,-0.2);
    \draw[thick](0,0) -- (-0.2,0.2);
  \fill[black] (1) circle (2pt) ;
  \fill[black] (2) circle (2pt) ; 
  \fill[black] (3) circle (2pt) ;
  \fill[black] (4) circle (2pt) ;  
  \draw[-stealth,thick] (2,0.5) -- (3, 0.5);
\coordinate (6) at (4, 0) ;
  \coordinate (7) at (5, 0) ;
  \coordinate (8) at (4, 1) ;
  \coordinate (9) at (5, 1) ;
  \coordinate (5) at (4.5, 1/2) ;
  \draw[line width=2.pt, red] (6) -- (5) ;
  \draw[line width=2.pt, red] (5) -- (7) ;
  \draw[line width=2.pt, red] (8) -- (5) ; 
  \draw[line width=2.pt, red] (5) -- (9) ; 
  \fill[black] (6) circle (2pt) ;
  \fill[black] (7) circle (2pt) ; 
  \fill[black] (8) circle (2pt) ;
  \fill[black] (9) circle (2pt) ;
  \fill[red] (5) circle (2pt) ; 
      \draw[thick](4,1) -- (4.2,1.2);
    \draw[thick](4,1) -- (3.8,1.2);
    \draw[thick](4,1) -- (3.8,0.8);
    \draw[thick](5,1) -- (4.8,1.2);
    \draw[thick](5,1) -- (5.2,1.2);
    \draw[thick](5,1) -- (5.2,0.8);
    \draw[thick](5,0) -- (5.2,0.2);
    \draw[thick](5,0) -- (5.2,-0.2);
    \draw[thick](5,0) -- (4.8,-0.2);
    \draw[thick](4,0) -- (4.2,-0.2);
    \draw[thick](4,0) -- (3.8,-0.2);
    \draw[thick](4,0) -- (3.8,0.2);
    \end{tikzpicture}
    \caption{The edge twist.}
    \label{fig:edge-twist}
\end{figure}

Generally, it is quite hard to describe the geometry of the resulting polyhedron and the arithmeticity of the associated reflection group.  
However, in some situations we are even able to decide (quasi-)arithmeticity of twisted antiprisms.  Let $A_{n,k}$, where $n \geqslant 4$ and $3 \leqslant k \leqslant n/2 + 1$, be an ideal hyperbolic right-angled polyhedron obtained from the antiprism $A_n$ by twisting two edges in one of its $n$-gonal faces such that the there are $k-2$ edges between the two twisted ones. Examples of $A_{4,3}$, $A_{6,3}$ and $A_{6,4}$ are presented in Figures~\ref{fig:A4} and~\ref{fig:A6}.

\begin{theorem}\label{th:arith-twisted-an}
   Let $A_{n,k}$, where $n \geqslant 4$ and $3 \leqslant k \leqslant n/2 + 1$, be an ideal hyperbolic right-angled twisted antiprism. Then the following properties hold. 
    \begin{itemize}
        \item[\textnormal{(i)}] $A_{n,k}$ is glued from $A_k$ and $A_{n-k+2}$ along their common ideal triangular face.
        \item[\textnormal{(ii)}] The associated right-angled reflection group $\Gamma_{n,k}$ is arithmetic if and only if $n=6$, $k=4$, or $n=4$, $k=3$.
        \item[\textnormal{(iii)}] $\Gamma_{n,k}$ is properly quasi-arithmetic if and only if $n=10$, $k=6$.
    \end{itemize}
\end{theorem}
\begin{remark}
    One can keep applying edge twist operations to $A_{n,k}$ a few more times. Notice that the twisted antiprism $A_{n,k} = A_k \cup_\Delta A_{n-k+2}$ (where $\Delta$ denotes a common ideal triangular face of $A_k$ and $A_{n-k+2}$) has one ideal regular $k$-gonal face, one ideal regular $(n-k+2)$-gonal face, and one ideal regular $n$-gonal face (a base), while all other faces are ideal triangles. Performing edge twists in the first two faces, one consequently splits $A_{n,k}$ into the union of a larger number of antiprisms. However, if an edge twist operation was applied to another face of $A_{n,k}$, it may significantly change the geometry of the polyhedron without giving such a concrete description as we have in Theorem~\ref{th:arith-twisted-an}.
\end{remark}

\subsection*{Structure of the paper}
Section~\ref{sec:preliminary} contains preliminaries on convex hyperbolic polyhedra (Section \ref{sec:convex-polyhedra}), adjoint trace fields and ambient groups of hyperbolic lattices (Section \ref{sec:trace-fields}), Vinberg invariants of hyperbolic reflection groups (Section \ref{sec:reflection-groups}), and antiprisms and their twists (Section \ref{sec:antiprisms}). Section \ref{section:proofs} is devoted to the proofs of Theorems \ref{th:general-GPS}, \ref{th:arith-surfaces-incomm-blocks}, \ref{th:Coxeter-orbifolds}, \ref{th:arith-out-of-quasi}, and \ref{th:arith-twisted-an}. We discuss some examples, namely ideal hyperbolic polygons and their arithmetic properties, in Section~\ref{sec:examples}. Finally, in Section~\ref{sec:open-problems}, we leave the list of open problems.

\subsection*{Funding}
Our work is supported by the Theoretical Physics and Mathematics Advancement Foundation ``BASIS''. The work of A.V. was also partially supported by the Ministry of Science and Higher Education of Russia (agreement no. 075-02-2024-1437).

\subsection*{Acknowledgements}
We thank Sami Douba for helpful discussions, comments and remarks on the earlier version of our manuscript. We are grateful to Stepan Alexandrov for his help with figures. We thank Leone Slavich for fruitful discussions and Misha Belolipetsky for his comments on the paper. We are also grateful to Alan Reid for interesting discussions on open problems related to hyperbolic knots, links, and their quasi-arithmeticity. The Coxeter prism example in Theorem~\ref{th:arith-out-of-quasi} was discovered by Nic Brody.

\section{Preliminaries}\label{sec:preliminary}

\subsection{Convex hyperbolic polyhedra}\label{sec:convex-polyhedra} The hyperbolic $d$-space $\HH^d$ is known to be the unique simply connected complete Riemannian $d$-manifold with constant sectional curvature $-1$. It can be defined as follows. Let $\mathbb{R}^{d,1}$ be the real vector space $\mathbb{R}^{d+1}$ equipped with the standard Lorentzian scalar product of signature $(d,1)$, namely, 
$$(x,y)=-x_0 y_0 +x_1 y_1 + \dots + x_d y_d.$$

Let us consider the two-sheeted hyperboloid $\mathcal{H}=\{x \in \mathbb{R}^{d,1} \,|\, (x,x) = -1 \}$ with two connected components 
$$
\mathcal{H}^+ = \{x \in \mathcal{H} \,|\, x_0 > 0\} \text{ and } \mathcal{H}^- = \{x \in \mathcal{H} \,|\, x_0 < 0\}.
$$

Then one can define the $d$-dimensional hyperbolic space $\HH^d$ as the Riemannian manifold $\mathcal{H}^+$ with the metric $\rho$ induced by restricting the bilinear form $(x,y)$ to the tangent bundle $T \mathcal{H}^+$. This hyperbolic metric $\rho$ satisfies $\cosh \rho(x, y) = - (x, y)$.
{\em Hyperplanes} of $\HH^d$ are intersections of linear hyperplanes of $\R^{d,1}$ with $\mathcal{H}^+$, and are totally geodesic submanifolds of codimension $1$ in $\HH^d$. 

In the hyperboloid model, points from the ideal hyperbolic boundary correspond to isotropic vectors:
$$
\partial \HH^d = \{x \in  \R^{d,1} \mid (x,x)=0 \textrm{ and }\, x_{0} > 0\}/\R_+ \approx \mathbb{S}^{d-1}.
$$

A \emph{convex hyperbolic $d$-polyhedron} is the intersection, with non-empty interior, of a finite family of closed half-spaces in hyperbolic $d$-space $\HH^d$. A \emph{hyperbolic Coxeter $d$-polyhedron} is a convex hyperbolic $d$-polyhedron $P$ all of whose dihedral angles are integer sub-multiples of $\pi$, i.e. of the form $\pi/m$ for some integer $m \geqslant 2$. A generalized\footnote{A \textit{generalized convex polyhedron} $P$ is the intersection, with non-empty interior, of possibly infinitely many closed half-spaces in hyperbolic $d$-space such that every closed ball intersects only finitely many bounding hyperplanes of $P$} hyperbolic convex polyhedron is called \emph{right-angled} if all its dihedral angles are $\pi/2$. A  generalized convex polyhedron is said to be \textit{acute-angled} if all its dihedral angles do not exceed $\pi/2$. 

It is known that generalized Coxeter polyhedra are the natural fundamental domains of discrete groups genereted by reflections in spaces of constant curvature, see~\cite{Vin85}.

A convex $d$-polyhedron has \emph{finite volume} if and only if it is the convex hull of finitely many points of the closure $\overline{\HH^d} = \HH^d \cup \partial \HH^d$. A convex polyhedron is said to be \emph{ideal}, if all its vertices are {\em ideal}, i.e. belong to the ideal hyperbolic boundary $\partial \HH^d$. 

Two compact polytopes $P$ and $P'$ in Euclidean space $\E^d$ are \emph{combinatorially equivalent} if there is a bijection between their faces that preserves the inclusion relation. A combinatorial equivalence class is called a \emph{combinatorial polytope}. Note that if a hyperbolic polyhedron $P \subset \HH^d$ is of finite volume, then the closure $\overline{P}$ of $P$ in $\overline{\HH^d}$ is combinatorially equivalent to a compact polytope of $\E^d$. 

\medskip

The following theorem is a special, right-angled, case of Andreev's theorem (see \cite{Andreev2} and \cite{RHD07}).

\begin{theorem}\label{theorem:Andreev} 
Let $\mathscr{P}$ be a combinatorial $3$-polytope. There exists a finite--volume  right-angled hyperbolic $3$-polyhedron $P \subset \overline{\HH^3}$ that realizes $\mathscr{P}$ if and only if:
\begin{enumerate}
    \item $\mathscr{P}$ is neither a tetrahedron, nor a triangular prism;
    \item every vertex of $\mathscr{P}$ belongs to at most four faces;
    \item if $f$, $f'$, and $f''$ are faces of $\mathscr{P}$, and $e' = f \cap f'$, $e'' = f \cap f''$ are non-intersecting edges, then $f'$ and $f''$ do not intersect each other;
    \item there are no faces $f_1$, $f_2$, $f_3$, $f_4$ such that $e_i := f_i \cap f_{i + 1}$ (indices $\mathrm{mod}\, 4$) are pairwise non-intersecting edges of $\mathscr{P}$.
\end{enumerate}
\end{theorem}

In a right-angled polyhedron $P \subset \overline{\HH^3}$, a vertex $v$ lies in $\HH^3$ (i.e. is finite) if and only if it belongs to exactly three faces of $P$, and $v$ is ideal, i.e. $v \in \partial\HH^3$, if and only if it is contained in four faces of $P$. Thus, all vertices of ideal hyperbolic right-angled $3$-polyhedra are $4$-valent.

\subsection{Lattices, adjoint trace fields and ambient groups}\label{sec:trace-fields}

In this subsection we briefly review the definitions and facts we will need about lattices and their arithmetic properties, see \cite[Section 2]{BBKS} for a more detailed discussion.

Let $\mathrm{O}_{d,1} = \mathbf{O}(f, \R)$ be the orthogonal group of the form $f(x) = (x,x)$, and $\mathrm{O}'_{d,1} < \mathrm{O}_{d,1}$ be the subgroup (of index $2$) preserving $\mathcal{H}^+$. 
The group $\mathrm{O}'_{d,1}$ preserves the metric $\rho$ on $\HH^d$, and is in fact isomorphic to the full isometry group $\mathrm{Isom}(\mathbb{H}^d)$ of the latter. 

If $\Gamma < \mathrm{O}'_{d,1}$ is a lattice, i.e., if $\Gamma$ is a discrete subgroup of $\mathrm{O}'_{d,1}$ with a finite-volume fundamental polyhedron in $\HH^d$, then the quotient $M=\mathbb{H}^d/\Gamma$ is a complete finite-volume {\itshape hyperbolic orbifold}. If $\Gamma$ is torsion-free, then $M$ is a complete finite-volume Riemannian manifold, and is called a {\itshape hyperbolic manifold}.

\begin{remark}
    Throughout this paper, we identify $\Isom(\HH^d)$ with the adjoint algebraic group $\PO_{d,1}(\R) = \mathrm{PO}_{d,1}$ instead of the Lie group $\mathrm{O}'_{d,1}$ since the latter is not an algebraic group. Besides, in this paper we apply our tools to reflection groups, and these are easier be to viewed through the hyperboloid model, which gives a natural identification of $\Isom(\HH^d)$ with $\mathrm{PO}_{d,1}$.
\end{remark}

Let us now provide the definitions of different arithmeticity types of hyperbolic lattices. Suppose that $G$ is a connected, semisimple real Lie group, and $\Gamma < G$ a subgroup. Moreover, we assume that it is an algebraic group, i.e. $G = \mathbf{G}(\R)$ for a real algebraic group $\mathbf{G}$. The group $\Gamma$ is called an \emph{arithmetic subgroup} of $G$ if there exist an algebraic $\mathbb{Q}$-group $\mathbf{H}$ and a surjective morphism $\varphi:\mathbf{H}(\mathbb{R})^{\circ}\rightarrow G$ with compact kernel such that $\Gamma$ is commensurable with $\varphi(\mathbf{H}(\mathbb{Z}))$.  If $G$ is not connected and thus has finitely many connected components, then we call $\Gamma < G$ arithmetic if $\Gamma \cap G^{\circ}$ is an arithmetic subgroup of $G^{\circ}$.
By the theorem of Borel and Harish-Chandra \cite{BHC62},  we have that any arithmetic subgroup $\Gamma$ of a semisimple Lie group $G$ is a lattice in $G$.

On the other hand, a fundamental theorem of Vinberg \cite[Theorem~1]{Vin71} implies that if $\Gamma$ is a Zariski-dense discrete subgroup of a semisimple adjoint real algebraic group $G$ (and by Borel's density theorem, any lattice is known to be Zariski-dense in the identity component $G^\circ$), then: 
\begin{itemize}
    \item[(1)] The {\em adjoint trace field} $k = \Q(\{\tr(\mathrm{Ad}\,\gamma) \mid \gamma \in \Gamma\})$, where $$\mathrm{Ad}:\mathrm{G}\rightarrow \mathrm{GL}(\mathfrak{g})$$ is the adjoint action of $G$ on its Lie algebra, is an invariant of the commensurability class of $\Gamma$;
    \item[(2)] There exists a maximal $k$-defined algebraic group $\mathbf{G}$ such that $$\mathbf{G}(\mathbb{R})^o\\ \cong G^o,\ \mathbf{G}(\mathbb{R}) < G,\ \text{and}\ \Gamma < \G(k),\ \text{up to conjugation.}$$ 
    It is called the {\itshape ambient group} of $\Gamma$; 
    \item[(3)] The ambient group $\G$ is uniquely determined up to $k$-isomorphism by the commensurability class of $\Gamma$.
\end{itemize}

\begin{remark}\label{admissible-def}
We will now assume that $G$ is an absolutely simple, adjoint, real algebraic group.  This is true, for instance, for all $\mathrm{PO}_{d,1}$ with $d \geqslant 2$ and $d \neq 3$. Let $\Gamma$ be an arithmetic lattice in $G$ with adjoint trace field $k$ and ambient group $\G$.
Then the following holds:
\begin{itemize}
    \item[(1)] the adjoint trace field $k$ is totally real,
    \item[(2)] the corresponding ambient group $\mathbf{G}$ is {\em admissible} for $G$, i.e. $\G(\R)^o$ is isomorphic to $G^o$ and $\G^\sigma(\R)$ is a compact group for any non-identity embedding $\sigma \colon k \hookrightarrow \mathbb{R}$,
    \item[(3)] and $\Gamma$ is  commensurable with the image in $G$ of $\G(\mathcal{O}_k)$, where $\mathcal{O}_k$ is the ring of integers of the field $k$.
\end{itemize}
\end{remark}  

\begin{remark}
Moreover, Tits' classification \cite{Tits} of semisimple algebraic $k$-groups implies that there exist only three types (in the terminology of \cite{BBKS}) of arithmetic lattices in $\mathrm{PO}_{d,1}$:
\begin{itemize}
    \item[(1)] \textit{type-I} arithmetic lattices arising from admissible quadratic forms (these lattices exist for all $d \geqslant 2$),
    \item[(2)] \textit{type-II} arithmetic lattices that come from skew-Hermitian forms over quaternion algebras (type-II lattices exist in all odd-dimensional hyperbolic spaces: $d = 3, 5, 7, \ldots$), and 
    \item[(3)] \textit{type-III} arithmetic lattices that comprise one exceptional family in $\mathrm{PO}_{3,1}$ and another in $\mathrm{PO}_{7,1}$.
\end{itemize} 
The type-I and type-II arithmetic lattices satisfy the definition through admissible ambient groups given in Remark~\ref{admissible-def}.     
The algebraic group $\mathrm{PO}_{3,1}$ is not absolutely simple, and the exceptional family of type-III arithmetic lattices gives rise to lattices whose adjoint trace field is not totally real and ambient group is not admissible. 
\end{remark}

\begin{remark}
In low-dimensional geometry and topology, it is quite common to identify $\Isom^+(\HH^3)$ with the group $\mathrm{PSL}_2(\CC)$. For lattices in $\mathrm{PSL}_2(\CC)$, the Vinberg adjoint trace field from \cite{Vin71} and the {\em invariant trace field} defined by Reid \cite{Reid90} (i.e. the field $\Q(\mathrm{tr}\,\Gamma^{(2)})$ where $\Gamma^{(2)} = \langle \gamma^2 \mid \gamma \in \Gamma \rangle$) coincide, see Vinberg \cite{Vin93}. However, if one views $\Gamma$ as a lattice in $\mathrm{PO}_{3,1}(\R)$, then there is some difference between these two notions of trace fields. For example, any arithmetic lattice $\Gamma$ of type I or II in the Lie group $\mathrm{PSO}_{3,1}(\R) \cong \Isom^+(\HH^3)$ has totally real adjoint trace field $k$, while the adjoint trace field of the image of $\Gamma$ in $\mathrm{PSL}_2(\CC)$ is an imaginary quadratic extension of $k$, cf. \cite[Coro. 3.12]{BBKS} and \cite[Theorem 2.3]{ALR01}. The key reason for these phenomena is that $\mathrm{PO}_{3,1}(\R)$ is not an absolutely (almost) simple algebraic group. The reader may consult \cite[Section 3.3]{BBKS} for more information about the behaviour of commensurability invariants of arithmetic lattices in  $\Isom(\HH^3)$.
\end{remark}

\begin{remark}
In this paper we are only interested in those arithmetic lattices $\Gamma$ that contain reflections in hyperplanes. Such arithmetic lattices contain in fact infinitely many reflections, and, moreover, they contain Zariski-dense reflection subgroups (see \cite[Theorem 1.1]{BK24}). Then it follows, see \cite[Lemma 7]{Vin67}, that the ambient groups of such $\Gamma$ are obtained from quadratic forms, and thus have to be admissible. 
\end{remark}

Now we are ready to define quasi-arithmetic hyperbolic lattices that were introduced first by Vinberg \cite{Vin67}.

\begin{definition}
If a lattice $\Gamma$ is contained, up to conjugation, in $\G(k)$ for an admissible $k$-group $\G$ (see Remark~\ref{admissible-def}), then $\Gamma$ is said to be \textit{quasi-arithmetic}. 
\end{definition}

As mentioned in the introduction, for every $d \geqslant 2$, there exist in $\HH^d$
\begin{itemize}
\item arithmetic lattices,
\item {\em properly quasi-arithmetic} lattices, i.e. those that are quasi-arithmetic but not arithmetic,
\item and lattices that are not quasi-arithmetic.
\end{itemize}
In each of these three large classes of hyperbolic lattices we have infinitely many pairwise incommensurable groups.

\medskip

In the proof of Theorem \ref{th:general-GPS} in Section \ref{section:proofs}, we will also make use of some results concerning adjoint trace fields of totally geodesic hypersurfaces in hyperbolic manifolds. 
\begin{theorem}[Mila \cite{Mila-thesis}]\label{th:Mila}
    Let $M$ be a hyperbolic manifold with adjoint trace field $k$ and ambient $k$-group $\mathbf{G}$ such that $\mathbf{G}(\mathbb{R}) = \mathrm{PO}_{n,1}$ with $n \geqslant 3$. Let $N \subset M$
    totally geodesic subspace of codimension $1$, i.e. $N = H/\Gamma_0$ where $H$ is a hyperplane being a lift of $N$ in $\HH^n$. Denote the reflection with respect to $H$ by $r$. Then
\begin{enumerate}
    \item[\textnormal{(1)}] $r \in \mathbf{G}(k)$;
    \item[\textnormal{(2)}] there exists a quadratic form $f$ over $k$ such that $\mathbf{G} = \mathbf{PO}_f$;
    \item[\textnormal{(3)}] the adjoint trace field of $\Gamma_0$ is a subfield of $k$.
\end{enumerate}
\end{theorem}

\noindent Part (3) is not explicitly stated in \cite{Mila-thesis}, but follows from the argument therein, and we explain it below.

\begin{proof}
    We refer to Section 3.2 in the thesis \cite{Mila-thesis} of Mila. Part (1) of our theorem is proved in \cite[Proposition 3.9]{Mila-thesis}. Note that for arithmetic hyperbolic manifolds part (1) also follows from recent work Belolipetsky et al.~\cite{BBKS}. Part (2) is precisely the statement of \cite[Theorem 3.8]{Mila-thesis}. 
    
    Part (3) follows from the proof of \cite[Proposition 3.9]{Mila-thesis}. Indeed, consider the centralizer $\mathbf{Z}(r)$ of the reflection $r$ in $\mathbf{G}$. It is a $k$-group since $r$ is $k$-defined. In the notation of \cite[Proposition 3.9]{Mila-thesis}, the group $\mathbf{Z}(r)$ splits over $k$ as the direct product $\mathbf{G}_0 \times \langle r\rangle$, where $\mathbf{G}_0$ is the Zariski-closure of $\Gamma_0$ in $\mathbf{G}$. Indeed, $\mathbf{G}_0$ is a normal subgroup of index $2$ (because $\langle r\rangle = \mathbb{Z}/2\Z$) and therefore is $k$-defined. Thus, $\Gamma_0 < \mathbf{G}_0(k)$ which implies that the adjoint trace field of $\Gamma_0$ is a subfield of $k$.
\end{proof}

\begin{corollary}\label{cor:quasi-arithmetic-subspace}
    In Theorem \ref{th:Mila}, if the lattice $\Gamma$ is quasi-arithmetic over $k$, then the adjoint trace field of $\Gamma_0$ is precisely $k$.
\end{corollary}
\begin{proof}
    By \cite[Theorems 1.7 \& 1.8]{BBKS}, the totally geodesic sublattice $\Gamma_0$ is quasi-arithmetic, and by Theorem \ref{th:Mila}, part (3), the adjoint trace field of $\Gamma_0$ is some $k_0 \subset k$. Obviously, $\G_0 < \G$. Let us consider the non-trivial embedding $\sigma$ of $k$ such that $\sigma|_{k_0} = \mathrm{id}$. The quasi-arithmeticity of $\Gamma_0$ implies that $\G_0^{\sigma}(\R)$ is non-compact. On the other hand, we have that $\G_0^{\sigma}(\R)$ is contained in $\G^{\sigma}(\R)$ with the latter being compact (due to the quasi-arithmeticity of $\Gamma$), a contradiction.
\end{proof}

The next statement (Theorem~\ref{th:inf-many-QA}) is slightly aside of the main subject of this paper, but highlights the important difference in the behavior of properly quasi-arithmetic vs. arithmetic lattices. 

\begin{theorem}[Mila \cite{Mila-doc}]\label{th:inf-many-QA}
For each $d \ge 2$, totally real number field $k$ and admissible over $k$ quadratic form $f$ of signature $(d,1)$, there exists an infinite sequence $\{\Gamma_n\}$ of pairwise incommensurable properly quasi-arithmetic lattices $\Gamma_n < \Isom (\HH^d)$ sharing the same adjoint trace field $k$ and ambient group $\mathbf{PO}_f$.
\end{theorem}

\begin{remark}\label{rem:Mila}
    In fact, Mila \cite{Mila-doc} obtained a much stronger result for $d \ge 3$ in the noncompact case and for $d \ge 4$ in the compact case. Namely, these quasi-arithmetic lattices in Theorem~\ref{th:inf-many-QA} can be chosen so that they even have the same trace ring (another invariant introduced by Vinberg in the same paper \cite{Vin71}). We remark that the trace ring is a stronger commensurability invariant than the ambient group or adjoint trace field. Theorem~\ref{th:inf-many-QA} without the trace ring condition is indeed true for all $d \ge 2$ with the same line of argument as in Mila's paper \cite{Mila-doc}. 
\end{remark}

\subsection{Vinberg invariants of hyperbolic reflection groups}\label{sec:reflection-groups}
In this subsection we define commensurability invariants of Zariski-dense hyperbolic reflection groups. Let  $H_e = \{x \in \HH^d \mid (x,e)=0\}$ be a hyperplane in $\HH^d \subset \R^{d,1} $ whose linear span in $\R^{d,1}$ has normal vector $e \in \R^{d,1}$ with $(e,e)=1$, and $H_e^- = \{x \in \HH^d \mid (x,e) \le 0\}$ be the half-space associated with it. If
$P = \bigcap_{j=1}^N H_{e_j}^- $
is a finite-sided Coxeter polyhedron in $\HH^d$, 
then the matrix 
$$G(P) = \{g_{ij}\}^N_{i,j=1} = \{(e_i, e_j)\}^N_{i,j=1}$$ is its {\em Gram matrix}. We write $K(P) = \Q\left(\{g_{ij}\}^N_{i,j=1}\right)$ and denote by $k(P)$ the field generated by all possible cyclic products of the entries of $G(P)$. For convenience, the set of all cyclic products of entries of a given matrix $A = (a_{ij})^N_{i,j=1}$, i.e., the set of all possible products of the form $a_{i_1 i_2} a_{i_2 i_3} \ldots a_{i_k i_1}$, will be denoted by $\Cyc(A)$. Thus, we have $k(P) = \Q\left(\Cyc(G(P))\right) \subset K(P)$.

Let $\Gamma_P$ denote a group generated by reflections in all the walls of $P$. Next, for $\{i_1, \ldots, i_s\} \subset \{1, \ldots, N\}$, we set
$$
v_1 = 2e_1; \quad v_{i_1} = 2 g_{1i_1} e_{i_1}; \quad v_{i_1 \ldots i_s} = 2^s g_{1i_1} g_{i_1 i_2} \ldots g_{i_{s-1} i_s} e_{i_s}.
$$
Denote by $V$ the linear span of all $v_{i_1 \ldots i_s}$. Since $P$ is a non-degenerate polyhedron in $\R^{d,1}$, the dimension of $V$ is $d+1$. The quadratic form $Q_P$ defined as the restriction on $V$ of the quadratic form (defined on $\R^N$) associated with the matrix $G(P)$ will be referred to as the {\em Vinberg form} of $P$. The matrix of the form $Q_P$ can be easily obtained as the Gram matrix of those $v_{i_1 \ldots i_s}$ that form a basis of $V$.  One can check that $V$ is a $k(P)$-defined subspace of $\R^N$ and is $\Gamma_P$-invariant, and thus $Q_P$ is a $k(P)$-form invariant under the $\Gamma_P$-action.

The following statement follows from work of Vinberg~\cite[Section 4]{Vin71} and a recent result of Dotti \cite{dotti}. 

\begin{theorem}[Vinberg invariants]\label{Vin-field}
  Let $P$ be a Coxeter polyhedron in $\HH^d$ and suppose $\Gamma_P$ is the associated reflection group. If $\Gamma_P$ is Zariski-dense in $\Isom(\HH^d)$, then $k(P)$ is the adjoint trace field of $\Gamma_P$, and thus is its commensurability invariant. Moreover, the Vinberg form $Q_P$ gives rise to the Vinberg ambient algebraic group $\mathbf{G}_P = \mathbf{PO}_{Q_P}$, i.e. there is a basis in which $\Gamma_P < \mathbf{G}_P(k(P))$, and therefore $Q_P$, up to similarity over $k(P)$, is also a commensurability invariant of~$\Gamma_P$.
\end{theorem}
\begin{remark}
    For reflection groups $\Gamma_P$ in $\HH^3$, the connection between the Vinberg field $k = k(P)$ and the invariant trace field $\ell$ is very explicit. If $d$ is the discriminant of the Vinberg form $Q_P$, then the invariant trace field of the index-$2$ orientation preserving subgroup $\Gamma_P^+$, viewed as a subgroup of $\mathrm{PSL}_2(\CC)$, satisfies the identity $\ell = k(\sqrt{d})$, see \cite[Theorem 10.4.1]{MR03}.
\end{remark}

We will use (in the proof of Theorem~\ref{th:arith-twisted-an}; see Section~\ref{sec:proof-th:arith-twisted-an}) the following criterion for determining if two orthogonal groups are isomorphic (see \cite[2.6]{GPS87}). 

\begin{lemma}\label{isogeny}
  Let $k$ be a field of characteristic $0$, $m \ge 2$, and $Q, Q'$ be two non-degenerate quadratic forms on the space $k^m$. Then $\PO_Q$ is $k$-isogenous (or, equivalently, $k$-isomorphic, since both groups are adjoint) to $\PO_{Q'}$ if and only if $Q$ and $Q'$ are similar over $k$, that is, if and only if there is some $\lambda \in k^\times$ such that the forms $Q$ and $\lambda Q'$ are isometric over $k$. 
\end{lemma}

The following criterion allows us to determine if a given finite-covolume hyperbolic reflection group $\Gamma$ with fundamental chamber $P$ is arithmetic, quasi-arithmetic, or neither.

\begin{theorem}[Vinberg's arithmeticity criterion \cite{Vin67}]\label{V}
Let $\Gamma < \mathrm{O}'_{d,1}$ be a reflection group with finite-volume fundamental Coxeter polyhedron $P \subset \HH^d$. Then $\Gamma$ is arithmetic if and only if each of the following conditions holds:
\begin{itemize}
    \item[{\bf(V1)}] $K(P)$ is a totally real algebraic number field;
    \item[{\bf(V2)}] for any embedding $\sigma \colon K(P) \to \R$, such that 
    $\sigma\!\mid_{k(P)} \ne \id$, the matrix $G^\sigma(P)$ is positive semi-definite;
    \item[{\bf(V3)}] $\Cyc(2 \cdot G(P)) \subset \OOO_{k(P)}$.
\end{itemize}
The group $\Gamma$ is quasi-arithmetic if and only if it satisfies conditions {\em\textbf{(V1)}--\textbf{(V2)}}, but not necessarily {\em\textbf{(V3)}}.
\end{theorem}

\begin{remark}\label{rem:non-compact-QA}
    The proof of Vinberg's arithmeticity criterion shows that if a finite-volume Coxeter polyhedron $P$ is not compact, then quasi-arithmeticity of $\Gamma_P$ implies $k(P) = \Q$. 
\end{remark}

\subsection{Antiprisms and their twists}\label{sec:antiprisms}

\begin{figure}
    \centering
    \includegraphics[scale=0.35]{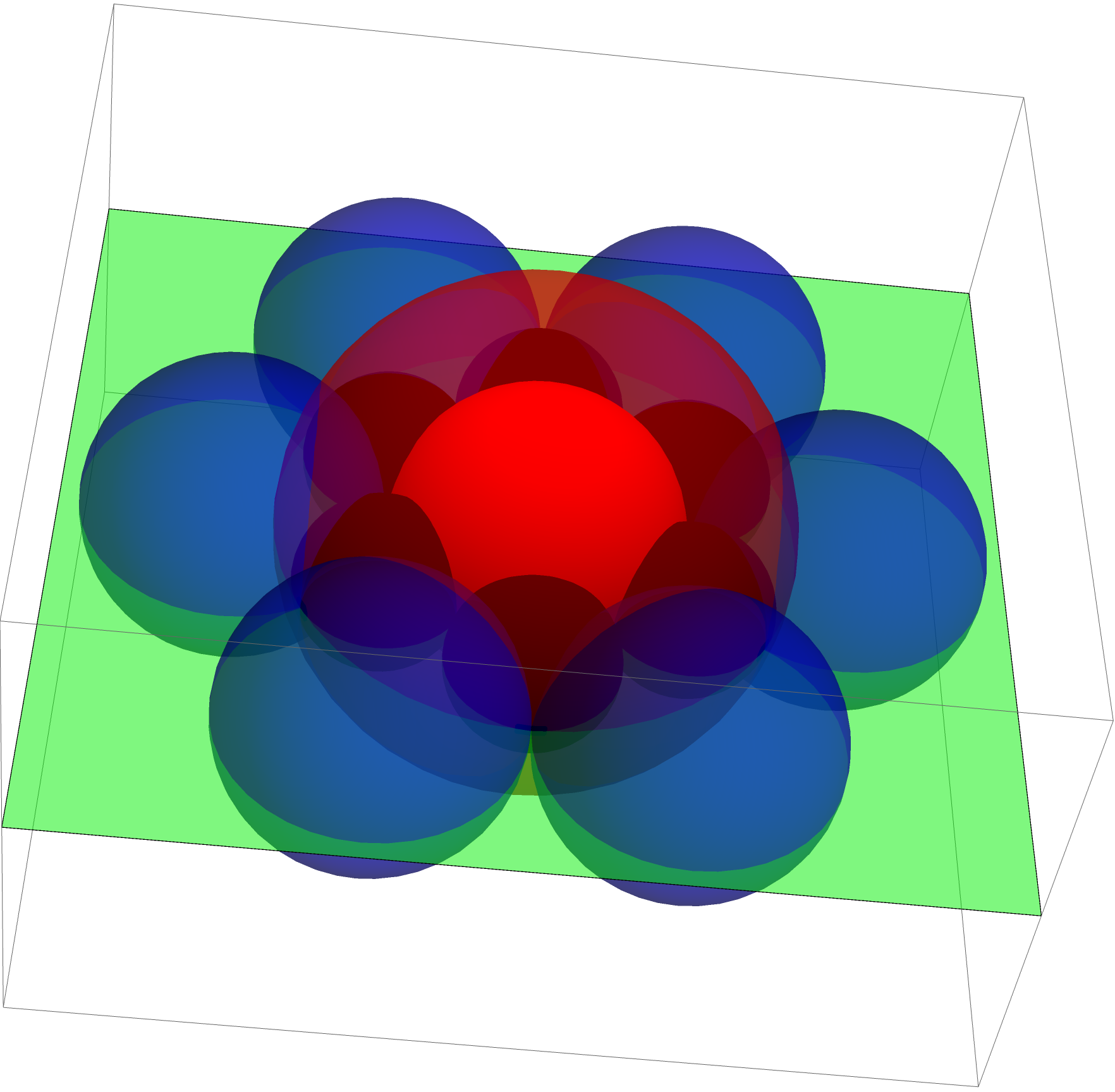}
    \caption{The antiprism $A_6$ in the upper half-space model of $\HH^3$ (where the green plane refers to the boundary) has two red ideal hexagonal bases, one is ``interior'' and another one ``exterior'', and twelve lateral faces bounded by blue spheres.}
    \label{fig:A6-3d}
\end{figure}

\begin{figure}
    \centering
    \includegraphics[scale = 0.3]{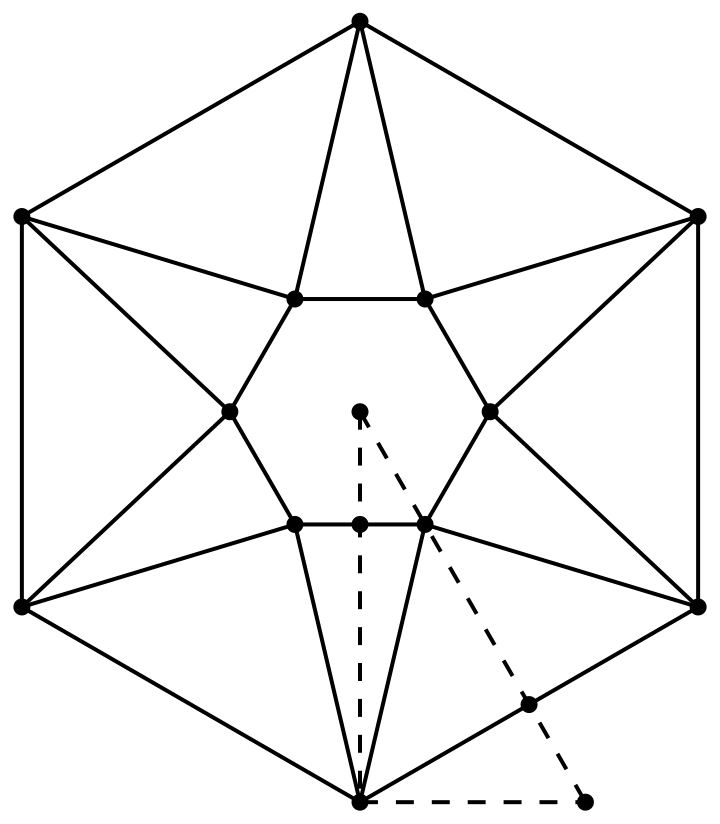} \quad
    \includegraphics[scale = 0.23]{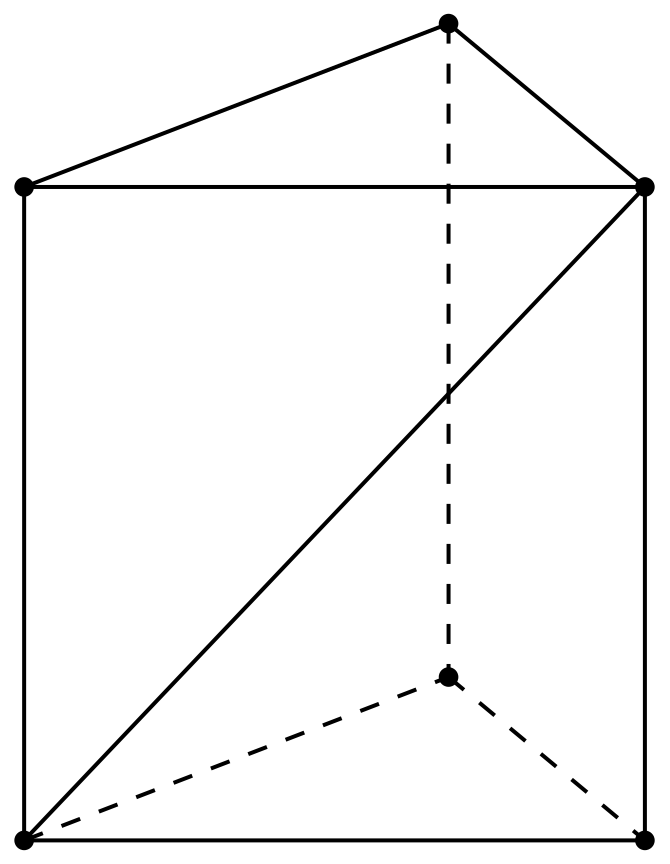}
    \caption{The antiprism $A_6$ and its slice $R_6$.}
    \label{fig:A6-cutting}
\end{figure}

The antiprism $A_n$, $n \geqslant 3$, is an ideal right-angled polyhedron in $\HH^3$ with $2n+2$ faces: two disjoint regular ideal $n$-gons, called {\em bases}, that are connected by an alternating band of $2n$ triangles. The visualization of the antiprism $A_6$ in the upper half-space model can be found in Figure~\ref{fig:A6-3d}, while its schematic picture is shown in Figure~\ref{fig:A6-cutting}. The schematic pictures of $A_3$ (an octahedron) and $A_4$ are depicted in Figure~\ref{fig:A4}. The realization of $A_n$ as an ideal right-angled polyhedron in $\HH^3$ can be easily verified by Andreev's theorem (see Theorem~\ref{theorem:Andreev}). To the best of our knowledge, the study of ideal right-angled antiprisms in the context of hyperbolic links was initiated by Thurston in~\cite[Ch.~6]{Th} where these polyhedra were called drums. Thurston remarked that the cases $n=3$ and $n=4$ are particularly interesting and important as the associated reflection groups of $A_3$ and $A_4$ are commensurable with $\operatorname{PSL}_2(\mathcal O_1)$ and $\operatorname{PSL}_2(\mathcal O_2)$, respectively. 

The polyhedron $A_n$ has many symmetries, since the centers of its $n$-gonal bases are the endpoints of the common perpendicular, and one base can be obtained from another via the $\frac{\pi}{n}$-twist around and shift along this perpendicular. This consideration allows to cut $A_n$ into $2n$ isometric pieces, each being a non-compact Coxeter polyhedron $R_n$ with $6$ faces; see Figure~\ref{fig:A6-cutting}. The commensurability of $A_n$ and $R_n$ is essential as it allows to explore arithmetic and group-theoretic properties of antiprisms. The key is that while it is complicated to provide the Gram matrix for the entire antiprism $A_n$, one can easily compute the Gram matrix of the $R_n$ (since the latter has the fixed amount of walls):
$$
\begin{pmatrix}
    1 & -\cos \frac{\pi}{n} & 0 & 0 & 0 & -1\\
    -\cos \frac{\pi}{n} & 1 & -1 & 0 & 0 & 0\\
    0 & -1 & 1 & -1 & 0 & 0\\
    0 & 0 & -1 & 1 & -\frac{1}{\cos \frac{\pi}{n}} & 0\\
    0 & 0 & 0 & -\frac{1}{\cos \frac{\pi}{n}} & 1 & -1\\
    -1 & 0 & 0 & 0 & -1 & 1
\end{pmatrix}
$$

This Gram matrix of $R_n$ is used to prove the following theorem.

\begin{theorem}\label{th:antiprisms}
    The following is true for antiprisms $A_n$:
\begin{enumerate}
    \item $k(A_n) = k(R_n) = \Q\left(\cos \frac{2\pi}{n}\right)$, $n \geqslant 3$;
    \item $A_n$ is arithmetic if and only if $n = 3,4$;
    \item $A_n$ is properly quasi-arithmetic antiprism if and only if $n = 6$;
    \item $A_n$ is not commensurable to $A_m$ for $n \neq m$.
\end{enumerate}
\end{theorem}

\noindent 
The parts $(1)$, $(2)$, and $(4)$ of Theorem \ref{th:antiprisms} are proved in \cite{Kel23} and \cite{MMT20}. The part~$(3)$ is not explicitly stated, but is in fact proved by Kellerhals in \cite{Kel23}. Indeed, if $A_n$ is quasi-arithmetic, then by Remark \ref{rem:non-compact-QA} we have $k(A_n) = \Q\left(\cos \frac{2\pi}{n}\right) = \Q$ which implies that $n = 3, 4, 6$. On the other hand, $A_n$ is quasi-arithmetic only if $R_n$ is so, and Kellerhals \cite{Kel23} shows that some of the cyclic products of the doubled Gram matrix of $R_6$ are rational but not integers. Thus, by Vinberg's arithmeticity criterion (Theorem~\ref{V}), we have that $A_6$ is properly quasi-arithmetic. We also note that the quasi-arithmeticity of the associated chain link $C_{12}$ was observed by Neumann--Reid \cite[Sections 7 \& 8]{NR90b}.

\begin{figure}[ht]
\begin{center}
\unitlength=.1mm
\begin{tikzpicture}[scale=0.38] 
\unitlength=10.mm
\draw [line width=2pt, black] (-16,6) circle[radius=0.1cm];
\draw [line width=2pt, black] (-15,4) circle[radius=0.1cm];
\draw [line width=2pt, black] (-17,4) circle[radius=0.1cm];
\draw [line width=2pt, black] (-16,1) circle[radius=0.1cm];
\draw [line width=2pt, black] (-20,7) circle[radius=0.1cm];
\draw [line width=2pt, black] (-12,7) circle[radius=0.1cm];
\draw [thick, black] (-20,7)-- (-12,7) -- (-16,1) -- (-20,7);
\draw [thick, black] (-15,4)-- (-16,6) -- (-17,4) -- (-15,4);
\draw [thick, black] (-16,1)-- (-17,4) -- (-20,7) -- (-16,6) -- (-12,7) -- (-15,4) -- (-16,1);
\draw(-16,-1) node {$A_3$}; 
\draw [line width=2pt, black] (0,5) circle[radius=0.1cm];
\draw [line width=2pt, black] (-8,5) circle[radius=0.1cm];
\draw [line width=2pt, black] (-4,1) circle[radius=0.1cm];
\draw [line width=2pt, black] (-4,9) circle[radius=0.1cm];
\draw [thick, black] (0,5)-- (-4,9) -- (-8,5) -- (-4,1) -- (0,5);
\draw [line width=2pt, black] (-3,4) circle[radius=0.1cm];
\draw [line width=2pt, black] (-5,4) circle[radius=0.1cm];
\draw [line width=2pt, black] (-3,6) circle[radius=0.1cm];
\draw [line width=2pt, black] (-5,6) circle[radius=0.1cm];
\draw [thick, black] (-3,4)-- (-3,6) -- (-5,6) -- (-5,4) -- (-3,4);
\draw [thick, black] (-3,4)-- (0,5) -- (-3,6) -- (-4,9) --(-5,6) -- (-8,5) -- (-5,4) -- (-4,1) -- (-3,4);
\draw(-4,-1) node {$A_4$};
\draw [dashed, ultra thick, red] (2,5) -- (14,5);
\draw [line width=2pt, black] (8,1) circle[radius=0.1cm];
\draw [line width=2pt, black] (8,9) circle[radius=0.1cm];
\draw [thick, black] (4,5)-- (8,9) -- (12,5) -- (8,1) -- (4,5);
\draw [line width=2pt, black] (7,4) circle[radius=0.1cm];
\draw [line width=2pt, black] (9,4) circle[radius=0.1cm];
\draw [line width=2pt, black] (7,6) circle[radius=0.1cm];
\draw [line width=2pt, black] (9,6) circle[radius=0.1cm];
\draw [thick, black] (7,4)-- (8,5) -- (7,6);
\draw [thick, black] (9,4)-- (8,5) -- (9,6);
\draw [thick, black]  (7,6) -- (9,6);
\draw [thick, black] (9,4) -- (7,4);
\draw [thick, black] (7,4)-- (4,5) -- (7,6) -- (8,9) --(9,6) -- (12,5) -- (9,4) -- (8,1) -- (7,4);
\draw(8,-1) node {$A_{4,3}$};
\draw [line width=2.pt, black] (4,5) circle[radius=0.1cm];
\draw [line width=2.pt, black] (12,5) circle[radius=0.1cm];
\draw [line width=2.pt, black] (8,5) circle[radius=0.1cm];
\end{tikzpicture}
\end{center}
\caption{Antiprisms $A_3$,  $A_4$, and the twisted antiprism $A_{4,3}$. The cutting of $A_{4,3}$ along the dashed red line shows the decomposition $A_{4,3} = A_3 \cup A_3$.} \label{fig:A4}
\end{figure}

Recall that in the introduction we defined an edge twist operation (see Figure~\ref{fig:edge-twist}). Remarkably, any ideal right-angled polyhedron can be obtained from one of the antiprisms by a series of such operations (see the recent papers of Erokhovets \cite{Er19} and Vesnin \cite{Ves17}). 

Since an edge twist is defined for a pair of disjoint edges, it can not be applied to the antiprism $A_3$. Figure~\ref{fig:A4} demonstrates an edge twist operation for $A_4$, where the resulting twisted antiprism $A_{4,3}$ is a union of two copies of $A_3$ glued along an ideal triangle. Thus, $\Gamma_{4,3}$ is an index two subgroup of the reflection group $\Lambda_3$ in the walls of $A_3$, and $\Gamma_{4,3}$ is arithmetic since the antiprism $A_3$ is arithmetic. 

\begin{figure}[ht]
\begin{center}
\unitlength=.1mm
\begin{tikzpicture}[scale=0.34] 
\unitlength=10.mm
\draw [line width=2pt, black] (-7,7) circle[radius=0.1cm];
\draw [line width=2pt, black] (-9,7) circle[radius=0.1cm];
\draw [line width=2pt, black] (-8,9) circle[radius=0.1cm];
\draw [line width=2pt, black] (-8,1) circle[radius=0.1cm];
\draw [line width=2pt, black] (-6,5) circle[radius=0.1cm];
\draw [line width=2pt, black] (-10,5) circle[radius=0.1cm];
\draw [line width=2pt, black] (-7,3) circle[radius=0.1cm];
\draw [line width=2pt, black] (-9,3) circle[radius=0.1cm];
\draw [line width=2pt, black] (-12,6.5) circle[radius=0.1cm];
\draw [line width=2pt, black] (-4,6.5) circle[radius=0.1cm];
\draw [line width=2pt, black] (-12,3.5) circle[radius=0.1cm];
\draw [line width=2pt, black] (-4,3.5) circle[radius=0.1cm];
\draw [thick, black] (-8,1)-- (-9,3) -- (-12,3.5) -- (-10,5) -- (-12,6.5)  -- (-9,7) -- (-8,9) -- (-7,7) -- (-4,6.5) -- (-6,5) -- (-4,3.5) -- (-7,3) -- (-8,1);
\draw [thick, black] (-8,1) -- (-12,3.5)-- (-12,6.5) -- (-8,9) -- (-4,6.5) -- (-4,3.5) -- (-8,1);
\draw [thick, black] (-9,3) -- (-10,5)-- (-9,7) -- (-7,7) -- (-6,5) -- (-7,3) -- (-9,3);
\draw(-8,-1) node {$A_{6}$};
\draw [dashed, ultra thick, red] (11,6.75) -- (5,6) -- (-1,6.75);
\draw [line width=2pt, black] (5,6) circle[radius=0.1cm];
\draw [line width=2pt, black] (4,7) circle[radius=0.1cm];
\draw [line width=2pt, black] (6,7) circle[radius=0.1cm];
\draw [line width=2pt, black] (5,9) circle[radius=0.1cm];
\draw [line width=2pt, black] (5,1) circle[radius=0.1cm];
\draw [line width=2pt, black] (3,5) circle[radius=0.1cm];
\draw [line width=2pt, black] (7,5) circle[radius=0.1cm];
\draw [line width=2pt, black] (4,3) circle[radius=0.1cm];
\draw [line width=2pt, black] (6,3) circle[radius=0.1cm];
\draw [line width=2pt, black] (9,6.5) circle[radius=0.1cm];
\draw [line width=2pt, black] (1,6.5) circle[radius=0.1cm];
\draw [line width=2pt, black] (9,3.5) circle[radius=0.1cm];
\draw [line width=2pt, black] (1,3.5) circle[radius=0.1cm];
\draw [thick, black] (5,1)-- (6,3) -- (9,3.5) -- (7,5) -- (9,6.5)  -- (6,7) -- (5,9) -- (4,7) -- (1,6.5) -- (3,5) -- (1,3.5) -- (4,3) -- (5,1);
\draw [thick, black] (5,1) -- (9,3.5)-- (9,6.5) -- (5,9) -- (1,6.5) -- (1,3.5) -- (5,1);
\draw [thick, black] (6,3) -- (7,5); \draw[thick, black] (6,7) -- (4,7);  \draw[thick, black] (3,5) -- (4,3) -- (6,3);
\draw [thick, black] (4,7) -- (5,6) -- (6,7);
\draw [thick, black] (3,5) -- (5,6) -- (7,5); 
\draw(5,-1) node {$A_{6,3}$};
\draw [dashed, ultra thick, red] (18,10) -- (18,5) -- (18,0);
\draw [line width=2pt, black] (18,5) circle[radius=0.1cm];
\draw [line width=2pt, black] (17,7) circle[radius=0.1cm];
\draw [line width=2pt, black] (19,7) circle[radius=0.1cm];
\draw [line width=2pt, black] (18,9) circle[radius=0.1cm];
\draw [line width=2pt, black] (18,1) circle[radius=0.1cm];
\draw [line width=2pt, black] (16,5) circle[radius=0.1cm];
\draw [line width=2pt, black] (20,5) circle[radius=0.1cm];
\draw [line width=2pt, black] (17,3) circle[radius=0.1cm];
\draw [line width=2pt, black] (19,3) circle[radius=0.1cm];
\draw [line width=2pt, black] (22,6.5) circle[radius=0.1cm];
\draw [line width=2pt, black] (14,6.5) circle[radius=0.1cm];
\draw [line width=2pt, black] (22,3.5) circle[radius=0.1cm];
\draw [line width=2pt, black] (14,3.5) circle[radius=0.1cm];
\draw [thick, black] (18,1)-- (19,3) -- (22,3.5) -- (20,5) -- (22,6.5)  -- (19,7) -- (18,9) -- (17,7) -- (14,6.5) -- (16,5) -- (14,3.5) -- (17,3) -- (18,1);
\draw [thick, black] (18,1) -- (22,3.5)-- (22,6.5) -- (18,9) -- (14,6.5) -- (14,3.5) -- (18,1);
\draw [thick, black] (19,3) -- (20,5) -- (19,7); \draw[thick, black]  (17,7) -- (16,5) -- (17,3); 
\draw [thick, black] (17,7) -- (18,5) -- (19,7);
\draw [thick, black] (17,3) -- (18,5) -- (19,3);
\draw(18,-1) node {$A_{6,4}$};
\end{tikzpicture}
\end{center}
\caption{The antiprism $A_6$ and the twisted antiprisms $A_{6,3}$ and $A_{6,4}$. The cuttings along red dashed lines decompose $A_{6,3}$ into $A_3 \cup A_5$ and $A_{6,4}$ into  $A_4 \cup A_4$ along ideal triangles.} \label{fig:A6}
\end{figure}

Figure~\ref{fig:A6} demonstrates two edge twist operations for $A_6$. The first one results in a twisted antiprism $A_{6,3}$ which is a union of antiprisms $A_3$ and $A_5$ glued along a common ideal triangle. Another one gives a twisted antiprism $A_{6,4}$ which is a union of two copies of $A_4$ glued together along an ideal triangle. Thus, $\Gamma_{6,4}$ is an index two subgroup of $\Gamma_4$ and $\Gamma_{6,4}$ is arithmetic as $\Gamma_4$ is arithmetic.

\section{Proofs}\label{section:proofs}

\subsection{Proof of Theorem~\ref{th:general-GPS}}
(1). The fact that the group $\Gamma = \langle \Delta_1, \Delta_2\rangle$ is a lattice was established by Gromov and Piatetski-Shapiro \cite[Section 2.10]{GPS87}, and their argument works in our situation without any arithmeticity assumptions on building blocks.

We recall that $\Gamma_i = N_i \rtimes \Delta_i$, where the reflection group $N_i = \langle \gamma r \gamma^{-1} \mid \gamma \in \Gamma_i \rangle$ is the normal closure of $r$ in $\Gamma_i$, $D_i$ is a fundamental domain for $N_i$ (bounded by some walls $\gamma H$), and $\Delta_i=\{\gamma \in \Gamma_i \mid \gamma(D_i)=D_i\}$. The next lemma is crucial for our results.

\vspace{-1.5em}
\begin{adjustwidth}{0.7cm}{0.7cm}
\begin{lemma}\label{lem:delta-r}
 $\Gamma_i = \langle \Delta_i, r\rangle$.
\end{lemma}
\begin{proof}
We claim that the symmetry group $\Delta_i$ transitively permutes the bounding walls of $D_i$. Indeed, we first observe that the entire collection of (pairwise disjoint) walls $\Omega_i = \{\gamma H \mid \gamma \in \Gamma_i\}$ is $\Gamma_i$-invariant and we assume that $D_i$ is bounded by some collection $\{\gamma_s H\}$ of walls including the hyperplane $H$ itself. Thus, any $\gamma D_i$ is just some fundamental chamber for $N_i$ (possibly coinciding with $D_i$), and there is a unique $n \in N_i$ such that $\gamma D_i = n D_i$, which gives a unique decomposition $\gamma = n \delta$, where $n \in N_i$, $\delta \in \Delta_i$.

We need to show that for any pair $\gamma_1 H$ and $\gamma_2 H$ of bounding walls of $D_i$ there is an element $\delta \in \Delta_i$ such that $\delta (\gamma_1 H) = \gamma_2 H$. Indeed, the isometry $\gamma = \gamma_2 \gamma_1^{-1}$ either takes $D_i$ to another chamber $D'_i$ adjacent to $D_i$ along the facet $\gamma_2 H$ ($\gamma_2 H$, the $\gamma$-image of $\gamma_1 H$, belongs to the $\gamma$-image of $D_i$) or leaves $D_i$ invariant (in this case, $\gamma \in \Delta_i$).

If $D_i$ is not invariant under $\gamma$, we have $D'_i = \gamma D_i = r' D_i$ for the reflection $r' = \gamma_2 r \gamma_2^{-1} \in N_i$. The above discussion shows that $\gamma = r' \delta$ for some $\delta \in \Delta_i$. Then $\delta  = (r')^{-1} \gamma = r' \gamma \in \Delta_i$ satisfies 
$$
\delta(\gamma_1 H) = r' \gamma (\gamma_1 H) = r' \gamma_2 \gamma_1^{-1} (\gamma_1 H) = r' (\gamma_2 H) = \gamma_2 H,
$$
since the reflection $r' \in N_i$ fixes $\gamma_2 H$ pointwise. 
Hence $\Gamma_i = \langle \Delta_i, r\rangle$.    
\end{proof}    
\end{adjustwidth}
\medskip

Lemma~\ref{lem:delta-r} allows us to justify Remark~\ref{rem:gluing} and to complete the proof of part (1). Indeed, the hyperplane $H$ projects to the connected, one-sided (since the entire $H$ is fixed by the reflection $r$), totally geodesic codimension-$1$ suborbifold $N=H/\Gamma_0$ of the finite-volume orbifold $M_i = \HH^n/\Gamma_i$.  

Now, consider the infinite-volume orbifold $M^+_i$, $i=1,2$,  obtained from the finite-volume orbifold $M_i = \HH^n/\Gamma_i$ by cutting off the totally geodesic finite-volume codimension-$1$ suborbifold $N = H/\Gamma_0$. Lemma~\ref{lem:delta-r} shows that the universal cover of $M_i^+$ is $D_i$ and therefore $M_i^+ = D_i/\Delta_i$. Moreover, $M_i^+$ has just one infinite end that can be cut along $N$ to obtain the finite-volume orbifold $M_i'$ with totally geodesic boundary $\partial M'_i = N$ (since adding the single reflection $r$ turns $M_i^+$ into $M_i$ due to Lemma~\ref{lem:delta-r}). In other words, each $M'_i$ is obtained from $M_i = \HH^n/\Gamma_i$ by forgetting the singularity of the totally geodesic codimension-$1$ subspace $N=H/\Gamma_0$, which is the singular locus of the reflection $r \in \Gamma_1 \cap \Gamma_2$. This implies that the gluing $M = M_1' \cup_N M_2'$ with the orbifold fundamental group $\pi_1^{orb} (M)= \Gamma = \Delta_1 *_{\Gamma_0} \Delta_2$ (this amalgamated free product appears due to the ping-pong lemma) has finite volume.

\begin{figure}
    \centering
    \includegraphics[width=0.9\linewidth]{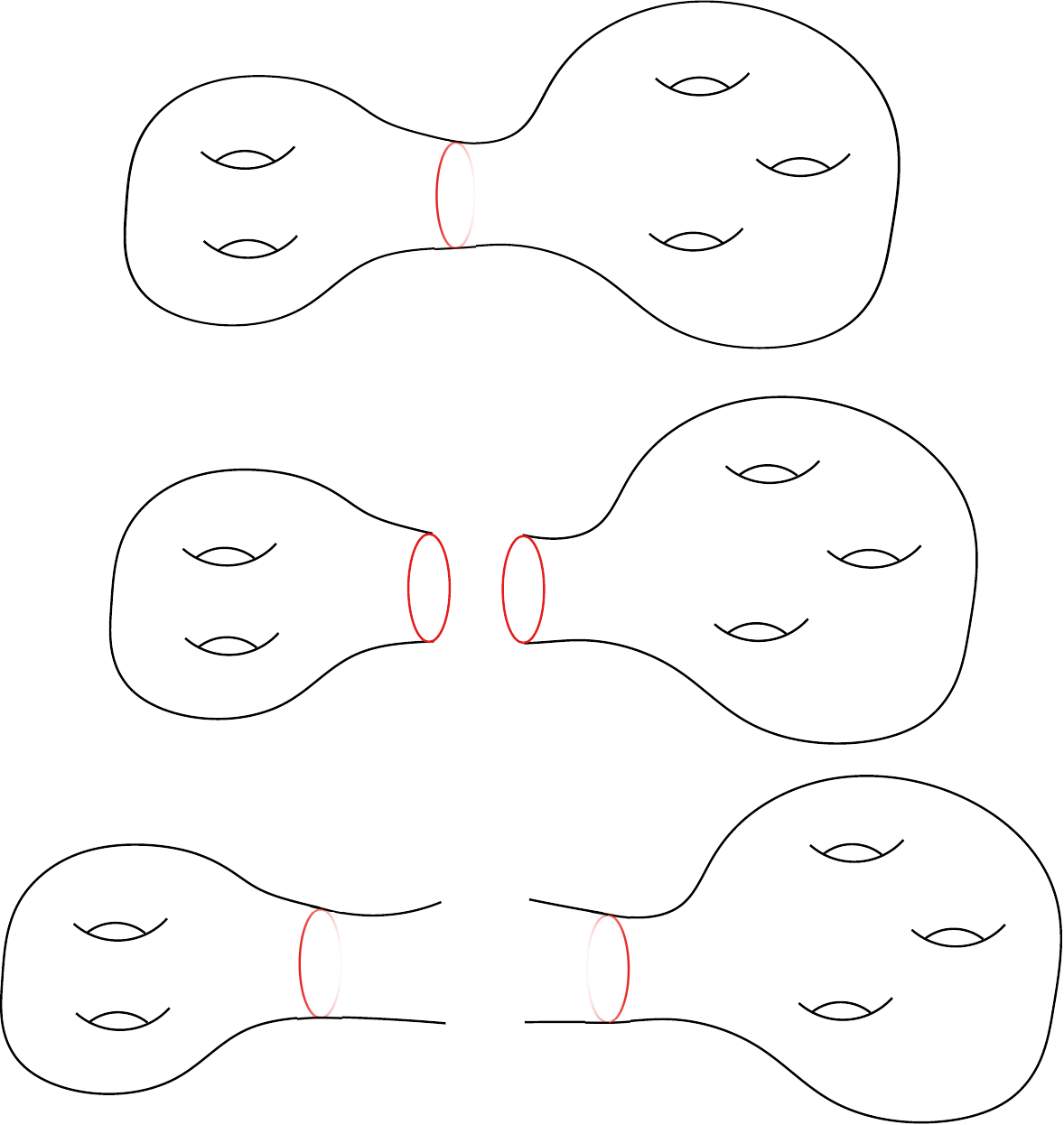}
    \caption{The gluing {\em manifold} $M$ is depicted on the top. The red arcs throughout this picture correspond to the totally geodesic hypersurface $N$. The shapes in the middle can be viewed as the building blocks, finite-volume manifolds $M_i'$ with isometric red geodesic boundaries, or as finite-volume orbifolds $M_i = \HH^n/\Gamma_i$ with the red curves being the singular loci (of $r$). Finally, infinite-volume manifolds $M_i^+=D_i/\Delta_i$ with infinite ends that can be cut along $N$ are shown in the bottom of this figure.}
    \label{fig:gluing-proof}
\end{figure}

In Figure~\ref{fig:gluing-proof}, we illustrate this situation for $M$ being a {\em manifold} of finite volume (in fact, one can always pass to a finite manifold cover by Selberg's lemma; in this case, $\Gamma$ and $\Delta_i$ become torsion-free). In the manifold case, the building blocks $M_i'$ are finite-volume {\em manifolds} with isometric red geodesic boundaries. As mentioned in Remark \ref{rem:gluing}, we conflate $M'_i$ with finite-volume {\em orbifolds} $M_i = \HH^n/\Gamma_i$ whose singular loci correspond to the totally geodesic submanifold $N$ and are the fixed point sets of the reflection $r \in \Gamma_1 \cap \Gamma_2$. Then $M_i^+=D_i/\Delta_i$ are infinite-volume manifolds with one infinite end each.

\medskip
(2). Now suppose that a gluing lattice $\Gamma$ is quasi-arithmetic with adjoint trace field $k$ and ambient group $\G$. The subgroups $\Delta_i < \Gamma$ are Zariski-dense, and thus have the same adjoint trace field $k$ and ambient group $\G$ (due to the admissibility of $\G$).

Then it remains to see that since $r \in \G(k)$ by Theorem~\ref{th:Mila}(1) and $\Delta_i < \G(k)$, we have (by Lemma~\ref{lem:delta-r}) $\Gamma_i = \langle \Delta_i, r\rangle < \G(k)$ and therefore $\Gamma_i$ is quasi-arithmetic with the same $k$ and $\G$. 

\medskip
(3). Since $\Gamma$ is nonarithmetic, the Margulis commensurator rigidity implies that $\mathrm{Comm}(\Gamma)$ is itself a lattice. Moreover, $\mathrm{Comm}(\Gamma)$ contains $r$ and each $\Delta_i$. Hence, the lattices $\Gamma_i = \langle \Delta_i, r \rangle$ (by Lemma~\ref{lem:delta-r}) are contained in the lattice $\mathrm{Comm}(\Gamma)$, whence they are finite-index subgroups of $\mathrm{Comm}(\Gamma)$. This implies that the lattices $\Gamma_1$ and $\Gamma_2$ are commensurable. \qed

\begin{remark}
Corollary \ref{cor:non-quasi-fields} easily follows from Theorem~\ref{th:general-GPS}, but we would also like to note that $H/\Gamma_0$ is a totally geodesic codimension-$1$ subspace of the orbifold $\HH^n/\Gamma$, and thus the lattice $\Gamma_0$ has the same adjoint trace field $k$ and ambient group $\G$ as $\Gamma$ by \cite[Theorems 1.7 \& 1.8]{BBKS} and Corollary \ref{cor:quasi-arithmetic-subspace}. Hence, one can rule out the conditions (1) and (2) in Corollary \ref{cor:non-quasi-fields} if the gluing $\Gamma$ is quasi-arithmetic.
\end{remark}

\subsection{Proof of Theorem~\ref{th:arith-surfaces-incomm-blocks}}
If $S_{g,1}$ is an arithmetic hyperbolic surface and $\alpha$ is a simple closed geodesic on it, then the main result of \cite{BBKS} shows that $\alpha$ is an fc-subspace of $S_{g,1}$ in the terminology of \cite{BBKS}, which means that $r \in \Comm(\pi_1 S_{g,1})$. 

On the other hand, since the closed geodesic $\alpha$ is separating, we obtain the splitting of $S_{g,1}$ into two orbifolds $M_1$ and $M_2$ (or, more precisely, the surfaces $M_1'$ and $M_2'$ with isometric boundaries, as in the proof of Theorem~\ref{th:general-GPS}), one being compact and another one having one cusp. This shows that the corresponding lattices $\Gamma_1$ and $\Gamma_2$ are not commensurable. \qed

\subsection{Proof of Theorem~\ref{th:Coxeter-orbifolds}}

Theorem~\ref{th:Coxeter-orbifolds} is a special case of Theorem~\ref{th:general-GPS}: $\Delta_i$ is just the group generated by reflections in all the walls of $P_i$ excepting the supporting hyperplane $H$ of the facet $F$. Then it remains to see from the construction, that $\Gamma_0$ is generated by reflections in hyperplanes that are orthogonal to $H$. 

Here we explicitly use the fact that $F$ meets its adjacent facets at even angles. Indeed, let $H$ and the supporting hyperplane $H'_j$ of a facet $F'_j$ of $P_j$ form the same angle $\pi/2m$ for $j=1, 2$. Then the stabiliser of $H \cap H'_j$ in both $\Gamma_{P_1}$ and $\Gamma_{P_2}$ is the dihedral group of order $2m$. This implies that there exists a reflection $r''$ in $\Gamma_{P_1} \cap \Gamma_{P_2}$  with mirror $H''$ orthogonal to $H$. This reflection $r''$ leaves $H$ invariant, and thus $r'' \in \Gamma_0$. This means that the stabiliser $\Gamma_0$ contains all reflections whose mirrors (intersected with $H$) bound the facet $F$ inside the hyperplane $H$. On the other hand, it is clear that the fundamental domain of $\Gamma_0$ in $H \cong \HH^{n-1}$ contains $F$, and thus coincides with $F$. Hence, $\Gamma_0$ is indeed generated by reflections in hyperplanes orthogonal to $H$.  \qed

\subsection{Proof of Theorem~\ref{th:arith-out-of-quasi}}

First of all, we note that DeBlois \cite{DB10} provided an example of a compact hyperbolic $3$-manifold $M$ with totally geodesic boundary such that its double across the boundary is nonarithmetic, but the gluing of $M$ with a certain twist $\tau$ of $M$ along $\partial M$ gives an arithmetic hyperbolic $3$-orbifold $N = M \cup_\tau M$ (see \cite[Proposition 0.3]{DB10}). By Theorem~\ref{th:general-GPS}, this manifold $M$ is quasi-arithmetic.

\medskip
For an orbifold example, consider a compact Coxeter prism $P$ in $\HH^3$ defined by the following Coxeter--Vinberg diagram:

\begin{center}
\begin{tikzpicture}
\coordinate (1) at (0, {-sqrt(3)/4}) ;
\coordinate (2) at (0, {sqrt(3)/4}) ;
\coordinate (3) at (1, 0) ;
\coordinate (4) at (2, 0) ;
\coordinate (5) at (3, 0) ;

\draw (0.05, {sqrt(3)/4}) -- (0.05, {-sqrt(3)/4}) ;
\draw (-0.05, {sqrt(3)/4}) -- (-0.05, {-sqrt(3)/4}) ;
\draw (0, {-sqrt(3)/4-0.05}) -- (1,-0.05) ;
\draw (0, {-sqrt(3)/4 + 0.04}) -- (1,0.04) ;
  \draw (2) -- (3) ; 
  \draw (3) -- (4) ; 
  \draw[dashed] (4) -- node [above] {$\sqrt{\nicefrac{7}{6}}$} (5) ;

  \fill[black] (1) circle (2pt) ;
  \fill[black] (2) circle (2pt) ; 
  \fill[black] (3) circle (2pt) ;
  \fill[black] (4) circle (2pt) ;
  \fill[black] (5) circle (2pt) ;
  
  \draw (1) circle (2pt) ;
  \draw (2) circle (2pt) ; 
  \draw (3) circle (2pt) ;
  \draw (4) circle (2pt) ;
  \draw (5) circle (2pt) ; 
\end{tikzpicture}    
\end{center}
It is easy to see from the diagram that the adjoint trace field of $P$ is $\Q$. By Vinberg's arithmeticity criterion (Theorem~\ref{V}), the reflection group $\Gamma_P$ is not quasi-arithmetic, but not arithmetic, since $\left(2\sqrt{\frac{7}{6}}\right)^2 = \frac{14}{3} \notin \Z$ (that is, the condition \textbf{(V3)} in Theorem~\ref{V} fails).

One can notice that this prism $P$ has a triangular base $F$ orthogonal to all its adjacent faces. The face has angles $(\pi/3, \pi/4, \pi/4)$ and thus has an order-$2$ symmetry allowing to glue $P$ with itself along $F$ in two different ways. One gluing is a standard reflective doubling of $P$ and thus also gives a quasi-arithmetic polyhedron. However, the second gluing gives already an arithmetic prism $P'$ (here we present a combinatorial picture of $P'$; the weights $m = 3, 4$ correspond to the dihedral angles $\pi/m$, and other angles are equal to $\pi/2$)

\begin{center}
\begin{tikzpicture}[scale=0.77]
\draw[thick, dashed] (-1, 0) -- (2, 0) ; 
\draw[thick] (-1, 0) -- (0.5, 2.6) ; 
\draw[thick, dashed] (2, 0) -- (0.5, 2.6) node[midway, right] {\large 3}; 

\draw[thick] (2, -1) -- (5, -1) node[midway, below] {\large 3}; 
\draw[thick] (2, -1) -- (3.5, 1.6); 
\draw[thick] (5, -1) -- (3.5, 1.6); 

\draw[thick] (-1, 0) -- (2, -1) node[midway, below] {\large 4}; 
\draw[thick, dashed] (2, 0) -- (5, -1) node[midway, above] {\large 3}; 
\draw[thick] (0.5, 2.6) -- (3.5, 1.6) node[midway, above] {\large 4}; 

\filldraw[black] (-1, 0) circle (2pt);
\filldraw[black] (2, 0) circle (2pt);
\filldraw[black] (0.5, 2.6) circle (2pt);
\filldraw[black] (2, -1) circle (2pt);
\filldraw[black] (5, -1) circle (2pt);
\filldraw[black] (3.5, 1.6) circle (2pt);

\end{tikzpicture}
\end{center}

with the Coxeter--Vinberg diagram

\begin{center}
\begin{tikzpicture}[scale=1.3]
\coordinate (1) at (1, {-sqrt(3)/4}) ;
\coordinate (2) at (1, {sqrt(3)/4}) ;
\coordinate (3) at (0, 0) ;
\coordinate (4) at (2, {-sqrt(3)/4}) ;
\coordinate (5) at (2, {sqrt(3)/4}) ;

\draw (1, {sqrt(3)/4-0.04}) -- (0,-0.04) ;
\draw (1, {sqrt(3)/4 + 0.06}) -- (0,0.06) ;
\draw (1, {-sqrt(3)/4-0.05}) -- (0,-0.05) ;
\draw (1, {-sqrt(3)/4 + 0.04}) -- (0,0.04) ;
  \draw (1) -- (2) ; 
  \draw (1) -- (4) ; 
  \draw (2) -- (5) ; 
  \draw[dashed] (4) -- node [right] {$\nicefrac{3}{2}$} (5) ;

  \fill[black] (1) circle (2pt) ;
  \fill[black] (2) circle (2pt) ; 
  \fill[black] (3) circle (2pt) ;
  \fill[black] (4) circle (2pt) ;
  \fill[black] (5) circle (2pt) ;
  
  \draw (1) circle (2pt) ;
  \draw (2) circle (2pt) ; 
  \draw (3) circle (2pt) ;
  \draw (4) circle (2pt) ;
  \draw (5) circle (2pt) ; 
\end{tikzpicture}    
\end{center}

Interestingly, this arithmetic prism $P'$ turns out to be a fundamental Coxeter polyhedron for the maximal reflection subgroup $R_f (\Z)$ preserving a quadratic form $f = \mathrm{diag}(-7,1,1,1)$ of signature $(3,1)$. Applying the program \texttt{VinAl} \cite{VinAl, BP18} (which is a software implementation of Vinberg's algorithm \cite{Vin72}) to the form $f$, one gets the following set of roots:

$e_1 = (0, -1, -1, 0)$, 

$e_2 = (0, 0, 1, -1)$,

$e_3 = (0, 1, 0, 0)$,

$e_4 = (1, -1, 2, 2)$,

$e_5 = (1, 0, 0, 3)$.

The corresponding Gram matrix is
$$
\begin{pmatrix}
2 & -1 & -1 & -1 & 0\\
-1 & 2 & 0 & 0 & -3\\
-1 & 0 & 1 & -1 & 0\\
-1 & 0 & -1 &  2 & -1\\
0 & -3 & 0 & -1 & 2
\end{pmatrix}
$$
One can easily check by reducing the norms of vectors $e_j$ that this Gram matrix corresponds to the Coxeter--Vinberg diagram of $P'$. \qed

\subsection{Proof of Theorem~\ref{th:arith-twisted-an}}\label{sec:proof-th:arith-twisted-an}

Let $e_1$ and $e_k$ be two disjoint edges of an $n$-gonal face of the antiprism $A_n$, $n \geqslant 4$, such that there are $(k-2)$ edges between $e_1$ and $e_k$. The combinatorial structure of $A_n$ gives us that the edge $e_1$ also belongs to some triangle, say $T_1$, of $A_n$. Denote by $V_1$ a vertex of $T_1$ opposite to $e_1$. Analogously, denote by $V_k$ a vertex opposite to $e_k$ in the triangle containing $e_k$. According to the definition of the edge twist operation, we firstly remove edges $e_1$ and $e_k$ and secondly, create a new vertex $V_0$ connected with the terminal vertices of $e_1$ and $e_k$. The resulting twisted antiprisms $A_{6,3}$ and $A_{6,4}$ are presented in Figure~\ref{fig:A6}. It can be seen from the pictures, that the ideal triangle $V_1 V_0 V_3$ decomposes the twisted antiprism $A_{6,3}$ into the antiprims $A_3$ and $A_5$, as well as the ideal triangle $V_1 V_0 V_4$ decomposes the twisted antiprism $A_{6,4}$ into two copies of the antiprism $A_4$. It is clear that for any $n$ and $k$ a common ideal triangle $V_1 V_0 V_k$ decomposes the twisted antiprism $A_{n,k}$ in antiprisms $A_k$ and $A_{n-k+2}$ which are ideal and right-angled by Theorem~\ref{theorem:Andreev}. This proves part (i).

To prove parts (ii) and (iii), we make use of Theorems~\ref{th:general-GPS} and \ref{th:Coxeter-orbifolds}. We first observe that the polyhedron $A_{n,k}$, as being glued from $A_k$ and $A_{n-k+2}$, has an ideal regular $k$-gonal face $F_1$ and an ideal regular $(n-k+2)$-gonal face $F_2$. If $\Gamma_{n,k}$ is quasi-arithmetic then both faces $F_1$ and $F_2$ must be quasi-arithmetic as well by \cite[Theorem 1.7]{BBKS} or \cite[Theorem 1.4]{BK21}. By Remark~\ref{rem:non-compact-QA}, the adjoint trace field of both $F_1$ and $F_2$ must be $\Q$, while for the regular ideal $M$-gon it is known to be isomorphic to $\Q(\cos \frac{2\pi}{M})$, since the reflection group in the sides of the regular ideal $M$-gon is commensurable with the triangle reflection group $(2, M, \infty)$. This leaves only the following options: $M = 3, 4, 6$. The latter implies that $k$ and $n-k+2$ can take only these values (recall that $k \leqslant n/2 + 1$ and therefore $k \leqslant n-k+2$):
$$
(k, n-k+2) = (3, 3), (3,4), (3,6), (4,4), (4,6), (6,6).
$$
Observe that the cases $(3,3)$ and $(4,4)$ correspond to arithmetic twisted antiprisms $A_{4,3} = A_3 \cup A_3$ and $A_{6,4} = A_4 \cup A_4$ being just the reflective doubles of the arithmetic antiprisms $A_3$ and $A_4$, respectively (their arithmeticity was established by \cite{MMT20, Kel23}; see Theorem~\ref{th:antiprisms}). The twisted antiprism $A_{10,6} = A_6 \cup A_6$ is commensurable to $A_6$ and the latter is properly quasi-arithmetic; see Theorem~\ref{th:antiprisms}, part (3). 

It remains to consider the twisted antiprisms 
$$
A_{5,3} = A_3 \cup A_4, \quad A_{7,3} = A_3 \cup A_6, \quad A_{8,4} = A_4 \cup A_6.
$$

The building blocks of these twisted antiprisms all have the same adjoint trace field $\Q$, and thus fall under conditions of Theorem~\ref{th:general-GPS}, part (2). We need to analyze the ambient groups $\mathbf{G}_3$, $\mathbf{G}_4$, and $\mathbf{G}_6$ of the antiprisms $A_3$, $A_4$, and $A_6$, respectively.

By Theorem~\ref{Vin-field}, it remains to compute the associated Vinberg forms. However, since the ambient group is a commensurability invariant, we can (and this is indeed easier) to pass to the slice $R_n$ of the antiprism $A_n$, see Section \ref{sec:antiprisms}. To that end, we  take the Gram matrix $G(R_n) = (g_{ij})$ of the slice $R_n$ and vectors
$$
v_1 = 2e_1, \quad v_2 = - 2 \cos \left(\frac{\pi}{n}\right) e_2,\quad v_{23} = 4\cos \left(\frac{\pi}{n}\right) e_3,\quad v_{234} = -8\cos \left(\frac{\pi}{n}\right) e_4.
$$
These vectors generate the Vinberg space $V$ for the reflection group of $R_n$. The Gram matrix of $R_n$ restricted onto $V$ gives a non-degenerate quadratic form
$$
Q_n = \begin{pmatrix}
    4 & 4\cos^2 \frac{\pi}{n} & 0 & 0\\
    4\cos^2 \frac{\pi}{n} & 4\cos^2 \frac{\pi}{n} & 8\cos^2 \frac{\pi}{n} & 0 \\
    0 & 8\cos^2 \frac{\pi}{n} & 16\cos^2 \frac{\pi}{n} & 32\cos^2 \frac{\pi}{n} \\
    0 & 0 & 32\cos^2 \frac{\pi}{n} & 64\cos^2 \frac{\pi}{n} \\
\end{pmatrix}
$$

Using the Gram--Schmidt orthogonalization process, we find the diagonal representative for each $Q_n$, where $n=3,4,6$:
$$
    Q_3 = \mathrm{diag}\left(-1, 1, 1, 1\right), \quad
    Q_4 = \mathrm{diag}\left(-2, 1, 1, 1\right), \quad
    Q_6 = \mathrm{diag}\left(-1, 1, 1, 3\right).
$$

The discriminants of these forms are $-1$, $-2$, $-3$, respectively. Now we will apply Lemma~\ref{isogeny}. Recall that if two quaternary forms $Q$ and $\lambda Q^\prime$ are similar then there is a transformation $C \in \mathrm{GL}_4(\Q)$ such that $\lambda Q^\prime = C^t \cdot Q \cdot C$ which implies
$$ 
\lambda^4 \det Q^\prime = (\det C)^2 \cdot \det Q.
$$
The latter means that
$$
\frac{\det Q^\prime}{\det Q} = \left(\frac{\det C}{\lambda^2}\right)^2 \in (\mathbb{Q}^*)^2.
$$

Observe that no one of the ratios of discriminants of the forms $Q_3$, $Q_4$, and $Q_6$ is a square of any rational number. Hence the forms $Q_3$, $Q_4$, and $Q_6$ are not pairwise similar, and the ambient groups $\mathbf{G}_3=\mathbf{PO}_{Q_3}$, $\mathbf{G}_4=\mathbf{PO}_{Q_4}$, and $\mathbf{G}_6=\mathbf{PO}_{Q_6}$ are not pairwise $\Q$-isogenous. Then Corollary \ref{cor:non-quasi-fields} from Theorem \ref{th:general-GPS} implies that the twisted antiprisms 
$$
A_{5,3} = A_3 \cup A_4, \quad A_{7,3} = A_3 \cup A_6, \quad A_{8,4} = A_4 \cup A_6
$$
are not quasi-arithmetic. \qed

\section{Arithmetic properties and gluings of ideal hyperbolic polygons}\label{sec:examples}

In this section, we discuss the arithmetic properties of groups generated by reflections in the sides of ideal hyperbolic polygons, as well as the behavior of gluings of such polygons. Recall that a convex polygon in $\HH^2$ is ideal if all its vertices belong to the boundary $\partial \HH^2$. We are going to use the upper half-plane model.

We also remark that all ideal triangles are isometric, and the corresponding ideal triangle reflection group $\Gamma_T$ is commensurable (in the wide sense) with $\mathrm{PSL}_2(\Z)$ and is therefore arithmetic. However, while any ideal $n$-gon $P$ in $\HH^2$ can be tesselated by $n-2$ ideal triangles, the corresponding reflection group $\Gamma_P$ is not necessarily commensurable with the ideal triangle reflection group $\Gamma_T$. In fact, $\Gamma_P$ may not even be quasi-arithmetic, since one can perform deformations of the distances between disjoint sides of $P$, which would affect the arithmeticity of $\Gamma_P$ due to Vinberg's arithmeticity criterion.

Nevertheless, there is a large family of ideal polygons whose reflection groups are quasi-arithmetic as we will see in the next subsection.

\subsection{Ideal polygons with rational vertices}

\begin{theorem}\label{th:rational-ideal}
    Let $r_1, \ldots, r_n \in \Q \cup \{\infty\}$ be pairwise distinct, and let $P$ be an ideal hyperbolic $n$-gon with vertices $r_1, \ldots, r_n$ in the upper half-plane model of $\HH^2$. Then the reflection group $\Gamma_P$ is quasi-arithmetic.
\end{theorem}

\noindent We will need a couple of formulas to prove this theorem.  We believe they are well-known, but we didn't find them in this specific form, so we provide justifications for the reader's convenience. These formulas will also be used in the next subsection.

\begin{lemma}\label{lem:distance-01x}
    Let $\ell_1$ be a vertical half-line $\mathrm{Re}(z)=0$, and $\ell_2$ be the hyperbolic geodesic with endpoints $1$ and $x \in \R_{>0}$ on the boundary line $\mathrm{Im}(z)=0$. 
\begin{center}
\begin{tikzpicture}
    \draw[thick,->] (-2,0) -- (3,0); 
    \fill[black] (0,0) circle (1.5pt) node [below] {$0$} ;
    \fill[black] (1,0) circle (1.5pt) node [below] {$1$} ;
    \fill[black] (2,0) circle (1.5pt) node [below] {$x$}; 
    \draw[thick,blue] (2,0) arc (0:180:0.5) node [midway, above] {{\small $\ell_2$}};

    \draw[thick, blue] (0,0) -- node [left] {{\small $\ell_1$}} (0,2);
\end{tikzpicture}
\end{center}
    Then the hyperbolic distance $\rho(\ell_1, \ell_2)$ between (disjoint) geodesics $\ell_1$ and $\ell_2$ satisfies
    $$
\cosh \rho(\ell_1, \ell_2) = \frac{x+1}{|x-1|} =
\begin{cases}
    1 + \frac{2}{|[0,\infty,1,x]|}, \quad \text{if } x>1, \\
    1 + \frac{2}{|[0,\infty,x,1]|}, \quad \text{if } x<1.
\end{cases}
    $$
\end{lemma}
\begin{proof}
    To find this distance, we will construct the common perpendicular $\ell$ to the lines $\ell_1$ and $\ell_2$. Since $\ell_1$ is a vertical Euclidean line, the geodesic $\ell$ is a circle centered at $0$ with radius $r$; see the picture below.
\begin{center}
\begin{tikzpicture}
    \draw[thick,->] (-3,0) -- (3,0); 

    \fill[black] (0,0) circle (1.5pt) node [below] {$0$} ;
    \fill[black] (1,0) circle (1.5pt) node [below] {$1$} ;
    \fill[black] (2,0) circle (1.5pt) node [below] {$x$}; 
    \fill[black] (4/3,{sqrt(2)/3}) circle (1.5pt) node [above right] {{\small $z = p+qi$}};
    \fill[black] (0,{sqrt(2)}) circle (1.5pt) node [above right] {{\small $r i$}};
    \draw[thick,red] ({sqrt(2)},0) arc (0:180:{sqrt(2)}) node [pos=0.3, above] {{\small $\ell$}};
    \draw[thick,blue] (2,0) arc (0:180:0.5) node [pos=0.05, above right] {{\small $\ell_2$}};

    \draw[thick, blue] (0,0) -- node [pos=0.9,left] {{\small $\ell_1$}} (0,2.5);
    \draw[thick] (0,0) -- node [above] {$r$} (4/3,{sqrt(2)/3});
    \draw[thick] (0,{sqrt(2)-0.2}) -- ++(0.2,0) -- ++(0,0.2);

    \begin{scope}[rotate around={atan(sqrt(2)/4):(4/3,{sqrt(2)/3})}]
        \draw[thick] (4/3-0.13,{sqrt(2)/3}) --  ++(0,0.15) -- ++(0.15,0);
    \end{scope}
\end{tikzpicture}
\end{center}

Calculating the degree of the point $0$ with respect to the circle $\ell_2$ shows that $r^2 = x$ (or one can consider similar Euclidean triangles $\Delta 01z$ and $\Delta 0zx$).
This implies that $r = \sqrt{x}$.

Then an elementary computation shows that the point of intersection of geodesics $\ell$ and $\ell_2$ has coordinates 
$$
p = \frac{2x}{x+1}, \quad q = \frac{\sqrt{x} \cdot |x-1|}{x+1}.
$$

Finally, using the well-known formula for the hyperbolic metric $\rho$ in the upper half plane, we obtain:
\begin{equation*}
\begin{split}
    \cosh \rho(\ell_1, \ell_2) = & \ \cosh \rho(ri, p+qi) = 1 + \frac{\|p+qi-ri\|^2}{2rq} \\
    = &\ 1+\frac{4x^2+((|x-1|-x-1)\sqrt{x})^2}{2x(x+1)|x-1|} = \frac{x+1}{|x-1|}.
\end{split}
\end{equation*}
It is easy to verify that $\frac{x+1}{|x-1|} = 1 + \frac{2}{x-1} = 1 + \frac{2}{|[0,\infty,1,x]|}$ if $x>1$. The case $x<1$ is similar.
\end{proof}

\begin{lemma}\label{lem:distance-geodesics}
 The following formulas for the hyperbolic distance $\rho(\ell_1, \ell_2)$ between disjoint geodesics $\ell_1$ and $\ell_2$ hold.
 
(1) Let $\ell_1$ be a vertical half-line $\mathrm{Re}(z)=a$, and $\ell_2$ be a hyperbolic geodesic with endpoints $b, c \in \R_{>0}$ on the boundary line $\mathrm{Im}(z)=0$, where $a < b < c$. Then 
    $$
\cosh \rho(\ell_1, \ell_2) = \frac{b+c-2a}{c-b} = 1 + \frac{2}{|[a,\infty,b,c]|}.
    $$
\begin{center}
\begin{tikzpicture}
    \draw[thick,->] (-2,0) -- (3,0); 
    \fill[black] (0,0) circle (1.5pt) node [below] {$a$} ;
    \fill[black] (1,0) circle (1.5pt) node [below] {$b$} ;
    \fill[black] (2,0) circle (1.5pt) node [below] {$c$}; 
    \draw[thick,blue] (2,0) arc (0:180:0.5) node [midway, above] {{\small $\ell_2$}};

    \draw[thick, blue] (0,0) -- node [left] {{\small $\ell_1$}} (0,2);
\end{tikzpicture}
\end{center}

(2) Let $\ell_1$ be a hyperbolic geodesic with endpoints $p, q \in \R$, $p<q$, and $\ell_2$ be a hyperbolic geodesic with endpoints $s, t \in \R$, $s<t$, on the boundary line $\R = \{z \mid \mathrm{Im}(z)=0\}$. Then 
    $$
\cosh \rho(\ell_1, \ell_2) = 
\begin{cases}
    1 + \frac{2}{|[p,q,t,s]|}, \quad \text{if } p<q<s<t\ \text{ (left picture)}, \\
    1 + \frac{2}{|[p,q,s,t]|}, \quad \text{if } p<s<t<q \ \text{ (right picture)}.
\end{cases}
    $$

\begin{center}
\begin{tikzpicture}
    \draw[thick,->] (-3,0) -- (3,0); 
    \fill[black] (-2,0) circle (2pt) node [below] {$p$} ;
    \fill[black] (0,0) circle (1.5pt) node [below] {$q$} ;
    \fill[black] (1,0) circle (1.5pt) node [below] {$s$} ;
    \fill[black] (2,0) circle (1.5pt) node [below] {$t$}; 
    \draw[thick,blue] (0,0) arc (0:180:1) node [midway, above] {{\small $\ell_1$}};
    \draw[thick,blue] (2,0) arc (0:180:0.5) node [midway, above] {{\small $\ell_2$}};

\end{tikzpicture}
\begin{tikzpicture}
    \draw[thick,->] (-3,0) -- (3,0); 

    \fill[black] (-1.5,0) circle (2pt) node [below] {$p$} ;
    \fill[black] (0,0) circle (1.5pt) node [below] {$s$} ;
    \fill[black] (1,0) circle (1.5pt) node [below] {$t$} ;
    \fill[black] (1.5,0) circle (1.5pt) node [below] {$q$}; 

    \draw[thick,blue] (1.5,0) arc (0:180:1.5) node [pos=0.75, above] {{\small $\ell_1$}};
    \draw[thick,blue] (1,0) arc (0:180:0.5) node [midway, above] {{\small $\ell_2$}};

\end{tikzpicture}
\end{center}
\end{lemma}
\begin{proof} In both cases, we find a Möbius transformation $F_i \in \mathrm{PSL}_2(\R)$ mapping four endpoints of the geodesics $\ell_1$ and $\ell_2$ to $0, \infty, 1, x$, and then we use Lemma~\ref{lem:distance-01x}.

(1) The Möbius transformation $F_1(z) = \frac{z-a}{b-a}$ sends $a \mapsto 0$, $\infty \mapsto \infty$, $b \mapsto 1$, and $c \mapsto x = \frac{c-a}{b-a}$. Note that $a<b<c$ implies that $x > 1$. Thus, with this $x$ we have $\cosh \rho(\ell_1, \ell_2) = \frac{x+1}{x-1} = \frac{b+c-2a}{c-b}$.

(2) The Möbius transformation
$$
F_2(z) = \frac{z-p}{z-q}\cdot \frac{s-q}{s-p}
$$
sends
$$
p \mapsto 0, \quad q \mapsto \infty, \quad s \mapsto 1, \quad t \mapsto x = \frac{(t-p)(s-q)}{(t-q)(s-p)}.
$$

Let $r_i$ denote the radius of the semicircle $\ell_i$. Then in the case $q<s<t$ (see the left side of the above picture for item (2)) we have $y := s-q > 0$ and
$$
x = \frac{(2r_1+2r_2+y)y}{(2r_1+y)(2r_2+y)} = \frac{y^2 + 2(r_1+r_2)y}{y^2 + 2(r_1+r_2)y + 4r_1 r_2} < 1.
$$
In the case $q>t>s$ (see the right side of the above picture for item (2)), we have $y := q-s > 0$ and
$$
x = \frac{(t-p)(q-s)}{(q-t)(s-p)} =\frac{(2r_1-y)(2r_2+y)}{y(2r_1-2r_2-y)} = \frac{-y^2 + 2(r_1-r_2)y + 4r_1 r_2}{-y^2 + 2(r_1-r_2)y } > 1.
$$
For the reader's convenience, we remark that in both cases the denominator and the numerator of the corresponding fractions for $x$ are positive.

Then it remains to note that the isometry $F_2$ preserves not only the distance, but also the cross-ratio of four points, and thus $[0,\infty,x,1]=[p,q,t,s]$ and similarly $[0,\infty,1,x]=[p,q,s,t]$.
\end{proof}

\medskip

\begin{proof}[Proof of Theorem~\ref{th:rational-ideal}]
    We will make use of Vinberg's arithmeticity criterion. Let $P$ be an ideal hyperbolic $n$-gon with vertices $r_1, \ldots, r_n \in \Q\cup \{\infty\}$ in the upper half-plane model of $\HH^2$. Denote the geodesic connecting $r_i$ and $r_{i+1}$ by $\ell_i$ for $i=1,\ldots,n$ (the geodesic $r_n r_1$ is denoted by $\ell_n$).

    Then the Gram matrix $G(P) = (g_{ij})^n_{i,j=1}$ of $P$ has the following entries:
    $$
g_{ii} = 1, \quad g_{i,i+1} = g_{i+1,i} = g_{1,n} = g_{n,1} = -1, \quad \text{ and all other } \ g_{ij} = -\cosh \rho(\ell_i, \ell_j). 
    $$
Lemma \ref{lem:distance-geodesics} shows that all $g_{ij} = -\cosh \rho(\ell_i, \ell_j)$ are rational functions of the endpoints $r_1, \ldots, r_n$. This implies $g_{ij} \in \Q$.

By Vinberg's arithmeticity criterion (Theorem~\ref{V}), the reflection group $\Gamma_P$ is quasi-arithmetic with adjoint trace field $\Q$.
\end{proof}

\begin{remark}
    Another approach to prove Theorem~\ref{th:rational-ideal} could be to pass to the orientation-preserving subgroup $\Gamma^+_P$ of the reflection group $\Gamma_P$. Then, one can use the fact that $\Gamma_P^+ < \mathrm{PSL}_2(\R)$ is generated by Möbius transformations that map geodesics with rational endpoints to geodesics with rational endpoints, because the index $[\Gamma_P:\Gamma^+_P] = 2$. This would imply that $\Gamma^+_P < \mathrm{PSL}_2(\Q)$ and thus is quasi-arithmetic. We provide another argument, which we believe is more rigorous and straightforward when one works with reflection groups.
\end{remark}

\subsection{Aritmetic ideal pentagon with a non-fc geodesic}

In this subsection, we provide an example of an arithmetic ideal pentagon $P$ glued from an arithmetic ideal triangle $T$ and a properly quasi-arithmetic ideal quadrilateral $P'$. That is, despite $\Gamma_P$ being arithmetic, the groups $\Gamma_1 = \Gamma_T$ and $\Gamma_2 = \Gamma_{P'}$ corresponding to the building blocks of $P$ are not commensurable. The key difference between this example and Theorem~\ref{th:arith-surfaces-incomm-blocks} is that $T$ and $P'$ are glued along an {\em infinite} common geodesic, which is fixed by the reflection $r$ that does not belong to $\Comm(\Gamma_P)$ (compare it with the proof of Theorem~\ref{th:arith-surfaces-incomm-blocks}).

Our construction is presented in the picture below. Consider an ideal pentagon $P$ bounded by five blue geodesics $\ell_1, \ldots, \ell_5$. We point out that geodesics $\ell_2$, $\ell_3$, and $\ell_4$ are semicircles of the same radius $a > 0$. This pentagon $P$ can be cut along the geodesic $\ell$ (the red arc) into the ideal triangle $T$ bounded by $\ell_1, \ell$, and $\ell_5$, and the ideal $4$-gon $P'$ bounded by $\ell$, $\ell_2$, $\ell_3$, and $\ell_4$.

\begin{center}
\begin{tikzpicture}
    \draw[thick,->] (-4,0) -- (4,0); 

    \fill[black] (0,0) circle (1.5pt) node [below] {$0$} ;
    \fill[black] (0.5,0) circle (1.5pt) node [below] {$a$} ;
    \fill[black] (1.5,0) circle (1.5pt) node [below] {$3a$}; 
    \fill[black] (-0.5,0) circle (1.5pt) node [below] {$-a$};
    \fill[black] (-1.5,0) circle (1.5pt) node [below] {$-3a$};
    \draw[thick,blue] (0.5,0)  arc (0:180:0.5) node [midway, below] {{\small $\ell_3$}} ;
    \draw[thick,blue] (-0.5,0) arc (0:180:0.5) node [midway, below] {{\small $\ell_2$}};
    \draw[thick,blue] (1.5,0) arc (0:180:0.5) node [midway, below] {{\small $\ell_4$}};
    \draw[thick,red]  (1.5,0) arc (0:180:1.5) node [midway, above] {{\small $\ell$}};

    \draw[thick, blue] (-1.5,0) -- node [left] {{\small $\ell_1$}} (-1.5,2.5);
    \draw[thick, dashed, blue] (-0.5,0) -- (-0.5,2.5);
    \draw[thick, dashed, blue] (0.5,0) -- (0.5,2.5);
    \draw[thick, blue] (1.5,0) -- node [right] {{\small $\ell_5$}}  (1.5,2.5);
\end{tikzpicture}
\end{center}

The blue dashed vertical half-lines at the points $\mathrm{Re}(z)=a$ and $\mathrm{Re}(z)=-a$ split $P$ into three ideal triangles, each of which can be mapped to the next by reflection with respect to one of these vertical dashed geodesics. This implies that $\Gamma_P$ is, in fact, commensurable to $\Gamma_T$, the ideal triangle reflection group, and thus is arithmetic.

On the other hand, $\rho(\ell,\ell_3) = \rho(ai, 3a i) = \log \left(\frac{3a}{a}\right)$ and thus $\cosh \rho(\ell,\ell_3) = \frac{5}{3}$ (we used the well-known formula $\rho(ai,bi) = \log(a/b)$, but one can also use Lemma~\ref{lem:distance-geodesics}). Also, by Lemma \ref{lem:distance-geodesics} we have 
$$
\cosh \rho(\ell_2, \ell_4) = 1 + \frac{2}{|[-3a,-a,3a,a]|} = 1+\frac{2\cdot 6a \cdot 2a}{2a \cdot 2a} = 7.
$$  
Applying Vinberg's arithmeticity criterion to the ideal quadrilateral $P'$, we see that $\Gamma_{P'}$ is quasi-arithmetic, but not arithmetic.

This means that $\Gamma_P$ and $\Gamma_T$ are not commensurable. We remark that the gluing locus of $T$ and $P$ is a bi-infinite red geodesic $\ell$ and thus it is not an fc-subspace of $P$ in the terminology of \cite{BBKS}. We recall that the results of Belolipetsky et al. \cite{BBKS} apply only to finite-volume totally geodesic subspaces of arithmetic orbifolds. In the case of $2$-dimensional orbifolds, this means that we can apply \cite{BBKS} to closed geodesics only (namely, in an arithmetic hyperbolic $2$-orbifold, every closed geodesic is an fc-subspace). In the next section, we ask if it is possible to construct a similar example with some of the blocks being even non-quasi-arithmetic.

\begin{remark}
    The arithmeticity of the ideal pentagon reflection group $\Gamma_P$ can also be verified directly by Vinberg's arithmeticity criterion. To use this criterion, it is enough to compute the distances $\cosh \rho(\ell_i, \ell_j)$. 
    
    We already know that $\cosh \rho(\ell_2, \ell_4) = 7$. Using Lemma \ref{lem:distance-geodesics}, we compute all the other distances:
    $$
\cosh \rho(\ell_1, \ell_3) = \cosh \rho(\ell_3, \ell_5) = 3, \quad
\cosh \rho(\ell_1, \ell_4) = \cosh \rho(\ell_2, \ell_5) = 5,
    $$
whence all entries of the Gram matrix $G(P)$ belong to $\Z$.
\end{remark}

\section{Open problems}\label{sec:open-problems}

Recent developments, including the results of this paper, highlight the importance of studying the arithmetic properties and commensurability invariants of hyperbolic lattices and manifolds. The first several open problems formulated below aim to attract attention to this general topic.

We recall that Vinberg \cite{Vin71} showed that the {\em (adjoint) trace ring} $A_\Gamma$ defined as the integral closure of $\Z[\{\tr\,\mathrm{Ad}(\gamma) | \gamma \in \Gamma \}]$ is an invariant of the  commensurability class of $\Gamma$. If the adjoint trace field $k=k_\Gamma$ is a number field, then we have
$$
A_\Gamma = \mathcal{O}_k[\{\tr\,\mathrm{Ad}(\gamma) | \gamma \in \Gamma \}]. 
$$
We would like to emphasize that the trace ring $A_\Gamma$ is a stronger commensurability invariant than $k_\Gamma$ or $\G_\Gamma$.
For example, any arithmetic lattice $\Gamma < \Isom(\HH^n)$ is completely determined up to commensurability by $k_\Gamma$ and $\G_\Gamma$, but even more can be said in terms of the trace rings: If for some other lattice $\Lambda < \Isom(\HH^n)$ we have $\G_\Lambda = \G_\Gamma$ and $A_\Lambda = A_\Gamma$, then $\Lambda$ is commensurable to the arithmetic lattice~$\Gamma$. 
Obviously, replacing the adjoint trace ring $A_\Gamma$ with the adjoint trace field $k_\Gamma$ in this statement leads to a false claim (in part due to the existence of properly quasi-arithmetic lattices). Another example of successful use of trace rings is due to Vinberg: in \cite{Vin12}, he constructed (via reflection groups) infinitely many properly quasi-arithmetic subgroups of $\mathrm{SL}_2(\Q)$. Vinberg showed that although these lattices have the same adjoint trace field and ambient group, their commensurability classes can be distinguished by adjoint trace rings.

\begin{question}
Do there exist any other commensurability invariant of hyperbolic lattices $\Gamma < \Isom(\HH^n)$, stronger than the adjoint trace ring $A_\Gamma$?   
\end{question}

\noindent We would like to stress that there is one quite strong commensurability invariant for {\em nonarithmetic} $1$-cusped hyperbolic orbifolds of finite volume. This invariant is called the {\em cusp density} (Neumann--Reid \cite[Proposition 1]{NR90} and Goodman--Heard--Hodgson \cite[Section 2]{GHH08}), and turned out to be very useful in certain situations. See, for example, Dotti--Drewitz--Kellerhals \cite{DDK23} and Kellerhals \cite{Kel23}.

\begin{question}
If a hyperbolic gluing $\Gamma$ in Theorem~\ref{th:general-GPS} is (quasi-)arithmetic, then do the corresponding (automatically quasi-arithmetic) lattices $\Gamma_1$ and $\Gamma_2$ have the same adjoint trace ring?
\end{question}

\noindent We refer to Mila \cite{Mila-doc} and Theorem \ref{th:inf-many-QA} (see also Remark~
\ref{rem:Mila}) to emphasize that there exist examples of pairwise incommensurable properly quasi-arithmetic lattices with the same adjoint trace ring. 

The next question is more specific and is related to our Theorem~\ref{th:arith-surfaces-incomm-blocks}.

\begin{question}
    Do there exist arithmetic hyperbolic gluings (compact ones in dimensions $\ge 2$ or finite-volume ones in dimensions $\ge 3$) with incommensurable building blocks? In other words, are there arithmetic hyperbolic manifolds with a separating totally geodesic hypersurface giving rise to incommensurable pieces (manifolds with boundary)?
\end{question}
\noindent In Theorem~\ref{th:arith-surfaces-incomm-blocks}, we provided $2$-dimensional non-compact such examples.

\begin{question}\textnormal{(A variation of Neumann--Reid \cite[Section 10, Question 3]{NR90b})}

\noindent Do there exist properly quasi-arithmetic hyperbolic knot complements?
\end{question}

\noindent As mentioned in the introduction, Reid \cite{Reid} proved that only one, the figure-eight-knot complement is arithmetic.

\begin{question}
    Is the number of commensurability classes of quasi-arithmetic hyperbolic link complements finite or infinite?
\end{question}

\noindent It is known that there are only finitely many commensurability classes of arithmetic hyperbolic link complements. More precisely, $14$ such classes are identified as follows from the work of Vogtman \cite{Vog85} on the cuspidal cohomology problem.

\begin{question}
    Is the number of commensurability classes of quasi-arithmetic ideal hyperbolic right-angled polyhedra in $\HH^n$, $n \ge 3$, finite or infinite?
\end{question}

\noindent The previous question is a special, right-angled, case of the more general open problem on finiteness of maximal quasi-arithmetic hyperbolic reflection groups; cf. \cite[Section 7, Question 5]{BK21}.

\begin{question}
    Extending the definition of quasi-arithmetic lattices to ambient groups of type II and III, we can ask the following. Do there exist general constructions of quasi-arithmetic lattices of types II and type III?
\end{question}

\noindent Alan Reid communicated to us that there is an example of a properly quasi-arithmetic lattice of type III in $\mathrm{PSL}_2(\CC)$. We expect to see more such examples among hyperbolic $3$-manifold groups. On the other hand, obtaining new families of such lattices in higher dimensions apparently requires an entirely new construction of nonarithmetic hyperbolic manifolds, which is a famous open problem both in real and complex hyperbolic geometry.

\bibliography{biblio.bib}{}

\begin{thebibliography}{10}

\bibitem{agol2006systoles}
{\sc I.~Agol}, {\em Systoles of hyperbolic 4-manifolds}, arXiv preprint
  math/0612290,  (2006).

\bibitem{ALR01}
{\sc I.~Agol, D.~D. Long, and A.~W. Reid}, {\em The {B}ianchi groups are
  separable on geometrically finite subgroups}, Ann. of Math. (2), 153 (2001),
  pp.~599--621.

\bibitem{All06}
{\sc D.~Allcock}, {\em Infinitely many hyperbolic {C}oxeter groups through
  dimension 19}, Geom. Topol., 10 (2006), pp.~737--758.

\bibitem{Andreev2}
{\sc E.~M. Andreev}, {\em Convex polyhedra of finite volume in
  {L}oba\v{c}evski\u{\i} space}, Mat. Sb. (N.S.), 83 (125) (1970),
  pp.~256--260.

\bibitem{AMR09}
{\sc O.~Antol\'{\i}n-Camarena, G.~R. Maloney, and R.~K.~W. Roeder}, {\em
  Computing arithmetic invariants for hyperbolic reflection groups}, in Complex
  dynamics, A K Peters, Wellesley, MA, 2009, pp.~597--631.

\bibitem{BBKS}
{\sc M.~Belolipetsky, N.~Bogachev, A.~Kolpakov, and L.~Slavich}, {\em Subspace
  stabilisers in hyperbolic lattices}, arXiv preprint arXiv:2105.06897,
  (2021).

\bibitem{BT11}
{\sc M.~V. Belolipetsky and S.~A. Thomson}, {\em Systoles of hyperbolic
  manifolds}, Algebr. Geom. Topol., 11 (2011), pp.~1455--1469.

\bibitem{BHW11}
{\sc N.~Bergeron, F.~Haglund, and D.~T. Wise}, {\em Hyperplane sections in
  arithmetic hyperbolic manifolds}, J. Lond. Math. Soc. (2), 83 (2011),
  pp.~431--448.

\bibitem{BD23}
{\sc N.~Bogachev and S.~Douba}, {\em Geometric and arithmetic properties of
  {L}\"obell orbifolds}, Algebr. Geom. Topol. (to appear), arXiv preprint
  arXiv:2304.12590,  (2023).

\bibitem{BDR24}
{\sc N.~Bogachev, S.~Douba, and J.~Raimbault}, {\em Infinitely many
  commensurability classes of compact coxeter polyhedra in $\mathbb{H}^4$ and
  $\mathbb{H}^5$}, Advances in Mathematics, 448 (2024), p.~109705.

\bibitem{BK21}
{\sc N.~Bogachev and A.~Kolpakov}, {\em On faces of quasi-arithmetic {C}oxeter
  polytopes}, Int. Math. Res. Not. IMRN,  (2021), pp.~3078--3096.

\bibitem{BK24}
\leavevmode\vrule height 2pt depth -1.6pt width 23pt, {\em Thin hyperbolic
  reflection groups}, arXiv preprint arXiv:2112.14642,  (2021), p.~10 pp.

\bibitem{VinAl}
{\sc N.~Bogachev and A.~Perepechko}, {\em Vinberg's algorithm}, Software
  implementation \texttt{VinAl},  (2017).

\bibitem{BP18}
{\sc N.~V. Bogachev and A.~Y. Perepechko}, {\em Vinberg's algorithm for
  hyperbolic lattices}, Mat. Zametki, 103 (2018), pp.~769--773.

\bibitem{BHC62}
{\sc A.~Borel and Harish-Chandra}, {\em Arithmetic subgroups of algebraic
  groups}, Ann. of Math. (2), 75 (1962), pp.~485--535.

\bibitem{CDBW12}
{\sc E.~Chesebro, J.~DeBlois, and H.~Wilton}, {\em Some virtually special
  hyperbolic 3-manifold groups}, Comment. Math. Helv., 87 (2012), pp.~727--787.

\bibitem{DB10}
{\sc J.~DeBlois}, {\em On the doubled tetrus}, Geom. Dedicata, 144 (2010),
  pp.~1--23.

\bibitem{dotti}
{\sc E.~Dotti}, {\em On the commensurability of hyperbolic {C}oxeter groups},
  Manuscripta Math., 173 (2024), pp.~203--222.

\bibitem{DDK23}
{\sc E.~Dotti, S.~T. Drewitz, and R.~Kellerhals}, {\em Cusp density and
  commensurability of non-arithmetic hyperbolic {C}oxeter orbifolds}, Discrete
  Comput. Geom., 69 (2023), pp.~873--895.

\bibitem{EM21}
{\sc V.~Emery and O.~Mila}, {\em Hyperbolic manifolds and
  pseudo-arithmeticity}, Trans. Amer. Math. Soc. Ser. B, 8 (2021),
  pp.~277--295.

\bibitem{Er19}
{\sc N.~Y. Erokhovets}, {\em Three-dimensional right-angled polytopes of finite
  volume in a {L}obachevski\u{\i} space: combinatorics and constructions}, Tr.
  Mat. Inst. Steklova, 305 (2019), pp.~86--147.

\bibitem{GHH08}
{\sc O.~Goodman, D.~Heard, and C.~Hodgson}, {\em Commensurators of cusped
  hyperbolic manifolds}, Experiment. Math., 17 (2008), pp.~283--306.

\bibitem{GPS87}
{\sc M.~Gromov and I.~Piatetski-Shapiro}, {\em Nonarithmetic groups in
  {L}obachevsky spaces}, Inst. Hautes \'{E}tudes Sci. Publ. Math.,  (1988),
  pp.~93--103.

\bibitem{Kel23}
{\sc R.~Kellerhals}, {\em A polyhedral approach to the arithmetic and geometry
  of hyperbolic link complements}, J. Knot Theory Ramifications, 32 (2023),
  pp.~Paper No. 2350052, 24.

\bibitem{Lack04}
{\sc M.~Lackenby}, {\em The volume of hyperbolic alternating link complements},
  Proc. London Math. Soc. (3), 88 (2004), pp.~204--224.
\newblock With an appendix by Ian Agol and Dylan Thurston.

\bibitem{MR03}
{\sc C.~Maclachlan and A.~W. Reid}, {\em The arithmetic of hyperbolic
  3-manifolds}, vol.~219 of Graduate Texts in Mathematics, Springer-Verlag, New
  York, 2003.

\bibitem{Mak68}
{\sc V.~S. Makarov}, {\em On {F}edorov's groups in four- and five-dimensional
  {L}obachevskij spaces}, Issled. po obshch. algebre. Kishinev (in Russian), 1
  (1968), pp.~120--129.

\bibitem{MMT20}
{\sc J.~S. Meyer, C.~Millichap, and R.~Trapp}, {\em Arithmeticity and hidden
  symmetries of fully augmented pretzel link complements}, New York J. Math.,
  26 (2020), pp.~149--183.

\bibitem{Mila18}
{\sc O.~Mila}, {\em Nonarithmetic hyperbolic manifolds and trace rings},
  Algebr. Geom. Topol., 18 (2018), pp.~4359--4373.

\bibitem{Mila-thesis}
\leavevmode\vrule height 2pt depth -1.6pt width 23pt, {\em The trace field of
  hyperbolic gluings}, PhD-Thesis, Universit{\"a}t Bern, 2019.

\bibitem{Mila-doc}
\leavevmode\vrule height 2pt depth -1.6pt width 23pt, {\em Non-commensurable
  hyperbolic manifolds with the same trace ring}, Doc. Math., 26 (2021),
  pp.~733--742.

\bibitem{Mila-imrn}
\leavevmode\vrule height 2pt depth -1.6pt width 23pt, {\em The trace field of
  hyperbolic gluings}, Int. Math. Res. Not. IMRN,  (2021), pp.~4392--4412.

\bibitem{NR90b}
{\sc W.~D. Neumann and A.~W. Reid}, {\em Arithmetic of hyperbolic manifolds},
  in Topology '90 ({C}olumbus, {OH}, 1990), vol.~1 of Ohio State Univ. Math.
  Res. Inst. Publ., de Gruyter, Berlin, 1992, pp.~273--310.

\bibitem{NR90}
\leavevmode\vrule height 2pt depth -1.6pt width 23pt, {\em Notes on {A}dams'
  small volume orbifolds}, in Topology '90 ({C}olumbus, {OH}, 1990), vol.~1 of
  Ohio State Univ. Math. Res. Inst. Publ., de Gruyter, Berlin, 1992,
  pp.~311--314.

\bibitem{Pur11}
{\sc J.~S. Purcell}, {\em An introduction to fully augmented links}, in
  Interactions between hyperbolic geometry, quantum topology and number theory,
  vol.~541 of Contemp. Math., Amer. Math. Soc., Providence, RI, 2011,
  pp.~205--220.

\bibitem{Pur20}
\leavevmode\vrule height 2pt depth -1.6pt width 23pt, {\em Hyperbolic knot
  theory}, vol.~209 of Graduate Studies in Mathematics, American Mathematical
  Society, Providence, RI, [2020] \copyright 2020.

\bibitem{Reid90}
{\sc A.~W. Reid}, {\em A note on trace-fields of {K}leinian groups}, Bull.
  London Math. Soc., 22 (1990), pp.~349--352.

\bibitem{Reid}
\leavevmode\vrule height 2pt depth -1.6pt width 23pt, {\em Arithmeticity of
  knot complements}, J. London Math. Soc. (2), 43 (1991), pp.~171--184.

\bibitem{RHD07}
{\sc R.~K.~W. Roeder, J.~H. Hubbard, and W.~D. Dunbar}, {\em Andreev's theorem
  on hyperbolic polyhedra}, Ann. Inst. Fourier (Grenoble), 57 (2007),
  pp.~825--882.

\bibitem{Ruz89}
{\sc O.~P. Ruzmanov}, {\em Examples of nonarithmetic crystallographic {C}oxeter
  groups in {$n$}-dimensional {L}obachevski\u{\i} space when {$6\leq n\leq
  10$}}, in Problems in group theory and in homological algebra ({R}ussian),
  Yaroslav. Gos. Univ., Yaroslavl, 1989, pp.~138--142.

\bibitem{Tho16}
{\sc S.~Thomson}, {\em Quasi-arithmeticity of lattices in
  {$\mathrm{PO}(n,1)$}}, Geom. Dedicata, 180 (2016), pp.~85--94.

\bibitem{Th}
{\sc W.~P. Thurston}, {\em The geometry and topology of three-manifolds. {V}ol.
  {IV}}, American Mathematical Society, Providence, RI, [2022] \copyright 2022.
\newblock Edited and with a preface by Steven P. Kerckhoff and a chapter by J.
  W. Milnor.

\bibitem{Tits}
{\sc J.~Tits}, {\em Classification of algebraic semisimple groups}, Algebraic
  Groups and Discontinuous Subgroups (Proc. Sympos. Pure Math., Boulder,
  Colo.), Amer. Math. Soc., Providence, RI,  (1966), pp.~33--62.

\bibitem{Ves91}
{\sc A.~Y. Vesnin}, {\em Three-dimensional hyperbolic manifolds with general
  fundamental polyhedron}, Math. Notes, 49 (1991), pp.~575--577.

\bibitem{Ves17}
{\sc A.~Y. Vesnin}, {\em Right-angled polytopes and three-dimensional
  hyperbolic manifolds}, Uspekhi Mat. Nauk, 72 (2017), pp.~147--190.

\bibitem{Vin67}
{\sc {\`E}.~B. Vinberg}, {\em Discrete groups generated by reflections in
  {L}oba\v{c}evski\u{\i}\ spaces}, Mat. Sb. (N.S.), 72 (114) (1967),
  pp.~471--488; correction, ibid. 73 (115) (1967), 303.

\bibitem{Vin71}
\leavevmode\vrule height 2pt depth -1.6pt width 23pt, {\em Rings of definition
  of dense subgroups of semisimple linear groups}, Math. of USSR--Izvestiya,
  {\bf 5} (1) (1971), pp.~45--55.

\bibitem{Vin72}
\leavevmode\vrule height 2pt depth -1.6pt width 23pt, {\em The groups of units
  of certain quadratic forms}, Mat. Sb. (N.S.), 87(129) (1972), pp.~18--36.

\bibitem{Vin85}
\leavevmode\vrule height 2pt depth -1.6pt width 23pt, {\em Hyperbolic
  reflection groups}, Russ. Math. Surv., 40 (1985), pp.~31--75.

\bibitem{Vin93}
\leavevmode\vrule height 2pt depth -1.6pt width 23pt, {\em The smallest field
  of definition of a subgroup of the group {${\rm PSL}_2$}}, Mat. Sb., 184
  (1993), pp.~53--66.

\bibitem{Vin12}
\leavevmode\vrule height 2pt depth -1.6pt width 23pt, {\em Some examples of
  {F}uchsian groups sitting in $\mathrm{SL}_2 (\mathbb{Q})$}, Preprint
  12–011, Universität Bielefeld: pp 4.
  http://www.math.uni-bielefeld.de/sfb701/files/preprints/sfb12011,  (2012).

\bibitem{Vin14}
\leavevmode\vrule height 2pt depth -1.6pt width 23pt, {\em Non-arithmetic
  hyperbolic reflection groups in higher dimensions}, Mosc. Math. J., 15
  (2015), pp.~593--602, 606.

\bibitem{Vog85}
{\sc K.~Vogtmann}, {\em Rational homology of {B}ianchi groups}, Math. Ann., 272
  (1985), pp.~399--419.

\end{thebibliography}
\bibliographystyle{siam}

\end{document}